\documentclass[10pt]{amsart}
\usepackage{latexsym, amsmath}

\renewcommand{\epsilon}{\varepsilon}
\newcommand{\newsection}[1]
{\subsection{#1}\setcounter{theorem}{0} \setcounter{equation}{0}
\par\noindent}

\newtheorem{theorem}{Theorem}

\newtheorem{lemma}[theorem]{Lemma}
\newtheorem{corr}[theorem]{Corollary}

\newtheorem{proposition}[theorem]{Proposition}
\newtheorem{deff}[theorem]{Definition}

\newcommand{\bth}{\begin{theorem}}
\newcommand{\ble}{\begin{lemma}}
\newcommand{\bcor}{\begin{corr}}

\newcommand{\bdeff}{\begin{deff}}

\newcommand{\bprop}{\begin{proposition}}
\newcommand{\eth}{\end{theorem}}
\newcommand{\ele}{\end{lemma}}
\newcommand{\ecor}{\end{corr}}
\newcommand{\edeff}{\end{deff}}

\newcommand{\eprop}{\end{proposition}}

\newcommand{\cd}{\, \cdot\, }

\renewcommand{\Pi}{\varPi}

\renewcommand{\epsilon}{\varepsilon}

\newcommand{\Rplus}{{\Bbb R}_+}

\newcommand{\R}{{\mathbb R}}
\newcommand{\ext}{{\R^3\backslash\mathcal{K}}}

\begin{document}

\title[Global existence for wave equations with multiple speeds]
{Global existence of solutions to multiple speed systems of
quasilinear wave equations in exterior domains}
\thanks{The first and third authors were supported in part by the NSF}

\author{Jason Metcalfe}
\address{School of Mathematics, Georgia Institute of Technology,
Atlanta, GA  30332-0160}
\author{Makoto Nakamura}
\address{Graduate School of Information Sciences, Tohoku
University, Sendai 980-8579, Japan}
\author{Christopher D. Sogge}
\address{Department of Mathematics,  Johns Hopkins University,
Baltimore, MD 21218}

\maketitle

\newsection{Introduction}
The goal of this paper is to prove global existence of solutions
to quadratic quasilinear Dirichlet-wave equations exterior to a
class of compact obstacles.  As in Metcalfe-Sogge \cite{MS}, the
main condition that we require for our class of obstacles is
exponential local energy decay (with a possible loss of
regularity).  Our result improves upon the earlier one of
Metcalfe-Sogge \cite{MS} by allowing a more general null condition
which only puts restrictions on the self-interaction of each wave
family.
The nonrelativistic system that we study serves as a
simplified model for the equations of elasticity.  In Minkowski
space, such equations were studied and shown to have global
solutions by Sideris-Tu \cite{Si3}, Agemi-Yokoyama \cite{AY}, and
Kubota-Yokoyama \cite{KY}.

The null condition we use here is the natural one for systems of
quasilinear wave equations with multiple speeds.  Following an
observation of John and Shatah, this null condition is equivalent
to the requirement that no plane wave solution of the system is
genuinely nonlinear (see John \cite{John}, p. 23 for the
single-speed case and Agemi and Yokoyama \cite{AY} for the
multi-speed case).
In order to allow the more general null condition, instead of just
exploring a coupling between a low order dispersive estimate and
higher order energy estimates as in Metcalfe-Sogge \cite{MS}, we
must first develop a low order energy estimate and couple this
with a low order pointwise estimate on the gradient and higher
order energy estimates.  Thus, our approach is a blend of the ones
using pointwise estimates based on fundamental solutions (see
e.g., \cite{KSS}, \cite{KSS3}, \cite{KY}, \cite{MS}, \cite{So2},
\cite{Y}) and ones using more refined $L^2$ energy
estimates for lower order terms (see, e.g. \cite{Hidano},
\cite{Si}, \cite{Si2}, \cite{Si3}).  As in the approach first
developed in \cite{KSS2} weighted space-time $L^2$ estimates for
lower order terms will also play a key role in our arguments.


We will be using an exterior domain analog of Klainerman's commuting
vector fields method \cite{knull}
.
Here, we have
to restrict to the collection of vector fields that are ``admissible''
for boundary value problems, $\{Z, L\}$, where $Z$ denotes the
generators of the spatial rotations and space-time translations
\begin{equation}\label{1.1}
Z=\{\partial_i, x_j\partial_k-x_k\partial_j, \: 0\le i\le 3, 1\le
j<k\le 3\}
\end{equation}
and $L$ is the scaling vector field
\begin{equation}\label{1.2}
L=t\partial_t + r\partial_r.
\end{equation}
Here and in what follows, $r=|x|$, and we will write
\begin{equation}\label{1.3}
\Omega_{ij}=x_i\partial_j-x_j\partial_i,\quad 1\le i<j\le 3.
\end{equation}
The generators of the hyperbolic rotations, $x_i\partial_t +
t\partial_i$, have an associated speed and the coefficients are
unbounded on the boundary of the obstacle, and thus, they do not seem
appropriate for the problem in question.

In Minkowski space, since $[(\partial_t^2-\Delta), Z]=0$ and
$[(\partial_t^2-\Delta),L]=2(\partial_t^2-\Delta)$, we see that $Z$
and $L$ preserve the equation $(\partial_t^2-\Delta)u=0$.  This is no
longer the case in the exterior domain since the Dirichlet boundary
condition is not preserved.  For the vector fields $Z$, since their
coefficients remain small in a neighborhood of our compact obstacle,
this is fairly easy to get around.  On the other hand, since the
coefficients of $L$ are large near the obstacle as $t\to\infty$, we
must stick to estimates that require relatively few of the scaling
vector field.

As in 
\cite{KSS2}, \cite{KSS3}, we will use weighted
$L^2_tL^2_x$ estimates where the weight is just a negative power of
$\langle x \rangle = \langle r\rangle = \sqrt{1+r^2}$.  These
estimates are useful for handling the lower order terms that arise in
the study of such boundary value problems.  They permit us to use
pointwise estimates for linear, inhomogeneous wave equations with
$O(\langle x\rangle^{-1})$ decay rather than the more standard
$O(t^{-1})$ decay which is more difficult to prove in the obstacle setting.
Additionally, such estimates allow us, as in 
\cite{MS}, to handle
the boundary terms that arise in the energy estimates if the obstacle
is no longer assumed to be star-shaped.  Here we exploit the fact that
we are studying equations with quadratic nonlinearities.

Additionally, we will be developing exterior domain analogs of a class
of weighted Sobolev estimates.  The weights here will involve powers
of $r$ and $\langle t-r\rangle$.  Specifically, we will be looking at
estimates of Klainerman-Sideris \cite{KS} and Hidano-Yokoyama
\cite{HY}.  We would additionally like to mention the works of Hidano
\cite{Hidano}, Kubota-Yokoyama \cite{KY}, Sideris \cite{Si, Si2},
and Sideris-Tu \cite{Si3} where similar estimates were used for the boundaryless case.

At this point, we wish to describe our assumptions on our obstacles
$\mathcal{K}\subset \R^3$.  We shall assume that $\mathcal{K}$ is
smooth and compact, but not necessarily connected.  By scaling,
without loss of generality, we may assume
\begin{equation*}
\mathcal{K}\subset \{\,x\in \R^3\,:\, |x|<1\}
\end{equation*}
throughout.  The only additional assumption states that there is
exponential local energy decay with a possible loss of regularity.
That is, if $u$ is a solution to $\Box u=0$ with Cauchy data $u(0,x)$,
$\partial_t u(0,x)$ supported in $|x|\le 4$, then there must be
constants $c, C>0$ so that
\begin{equation}\label{1.4}
\left(\int_{\{x\in\ext\,:\,|x|<4\}}
|u'(t,x)|^2\:dx\right)^{1/2}\le Ce^{-ct}\sum_{|\alpha|\le
1}\|\partial_x^\alpha u'(0,\cdot)\|_2.
\end{equation}
Here, and throughout, we are taking $\partial=\nabla_{t,x}$ to be the
space-time gradient.

We note that we do not require exponential decay; in fact,
$O((1+t)^{-1-\delta})$ may be sufficient with a tighter argument.  For
simplicity, we will assume \eqref{1.4}.  Currently, the authors are
not aware of any 3-dimensional example that involves polynomial decay,
but does not have exponential decay.

Notice that if the obstacle is assumed to be nontrapping, then a
stronger version of \eqref{1.4} holds where $\alpha=0$ on the right
side (see, e.g., Morawetz-Ralston-Strauss \cite{MRS}).  If there are
trapped rays, it was shown by Ralston \cite{ralston} that \eqref{1.4}
could not hold without a loss of regularity $\ell>0$ in the right
side.  We will assume throughout that $\ell=1$.  This can be done
without loss of generality since if $\ell>1$, interpolation with the
standard energy inequality will yield \eqref{1.4} (with a different
constant $c$).  In fact, we could take $\ell=\delta$ for any
$\delta>0$.

Ikawa \cite{Ikawa1}, \cite{Ikawa2} showed that there was such
exponential decay of local energy for certain finite unions of smooth,
convex obstacles with a loss $\ell=7$.  In particular, using Ikawa's
result, we have \eqref{1.4} for two disjoint convex obstacles or any
number of sufficiently separated balls.

For such smooth, compact obstacles $\mathcal{K}\subset\R^3$ satisfying
\eqref{1.4}, we shall consider quadratic, quasilinear
systems of the form
\begin{equation}\label{1.5}
\begin{cases}
\Box u = Q(du, d^2u),\quad (t,x)\in \Rplus\times\ext\\
u(t,\cd)|_{\partial\mathcal{K}}=0\\
u(0,\cd)=f,\quad \partial_t u(0,\cd)=g.
\end{cases}
\end{equation}
Here
\begin{equation}\label{1.6}
\Box=(\Box_{c_1},\Box_{c_2}, \dots, \Box_{c_D})
\end{equation}
is a vector-valued multiple speed d'Alembertian with
\begin{equation*}
\Box_{c_I}=\partial_t^2-c_I^2\Delta.
\end{equation*}
We will assume that the wave speeds $c_I$ are positive and
distinct.  This situation is referred to as the nonrelativistic
case.  Straightforward modifications of the argument give the more
general case where the various components are allowed to have the
same speed.  Also, $\Delta=\partial_1^2+\partial_2^2+\partial_3^2$ \ 
is the standard Laplacian.  Additionally, when convenient, we will
allow $x_0=t$ and $\partial_0=\partial_t$.

We shall assume that $Q(du,d^2u)$ is of the form
\begin{equation}\label{1.7}
Q^I(du,d^2u)=B^I(du)+\sum_{\substack{0\le j,k,l\le 3\\ 1\le J,K\le
D}}B^{IJ,jk}_{K,l}\partial_l u^K \partial_j\partial_k u^J, \quad
1\le I\le D
\end{equation}
where $B^I(du)$ is a quadratic form in the gradient of $u$ and
$B^{IJ,jk}_{K,l}$ are real constants satisfying the symmetry
conditions
\begin{equation}\label{1.8}
B^{IJ,jk}_{K,l}=B^{JI,jk}_{K,l}=B^{IJ,kj}_{K,l}.
\end{equation}
To obtain global existence, we shall also require that the
equations satisfy the following null condition which only involves
the self-interactions of each wave family.  That is, we require
that
\begin{equation}\label{1.9}
\sum_{0\le j,k,l\le 3}B^{JJ,jk}_{J,l}\xi_j\xi_k\xi_l=0
\quad\text{whenever}\quad
\frac{\xi_0^2}{c^2_J}-\xi_1^2-\xi_2^2-\xi_3^2=0,\quad J=1,\dots,D.
\end{equation}

To describe the null condition for the lower order terms, we expand
$$B^I(du)=\sum_{\substack{1\le J,K\le D\\ 0\le j,k\le
3}}A^{I,jk}_{JK}\partial_j u^J \partial_k u^K.$$ We then require
that each component satisfies the similar null condition
\begin{equation}\label{1.10}
\sum_{0\le j,k\le 3} A^{J,jk}_{JJ}\xi_j\xi_k=0 \quad\text{whenever }
\quad \frac{\xi_0^2}{c^2_J}-\xi_1^2-\xi_2^2-\xi_3^2=0, \quad J=1,\dots, D.
\end{equation}

In order to solve \eqref{1.5} we must also assume that the data
satisfies the relevant compatibility conditions.  Since these are
well-known (see, e.g., \cite{KSS}), we shall only describe them
briefly.  To do so we first let $J_ku=\{\partial^\alpha_x u : 0\le
|\alpha|\le k\}$ denote the collection of all spatial derivatives
of $u$ of order up to $k$.  Then if $m$ is fixed and if $u$ is a
formal $H^m$ solution of \eqref{1.5}, we can write
$\partial_t^ku(0,\cd)=\psi_k (J_k f, J_{k-1}g)$, $0\le k \le m$,
for certain compatibility functions $\psi_k$ which depend on the
nonlinear term $Q$ as well as $J_k f$ and $J_{k-1} g$.  Having
done this, the compatibility condition for \eqref{1.5} with
$(f,g)\in H^m\times H^{m-1}$ is just the requirement that the
$\psi_k$ vanish on $\partial\mathcal{K}$ when $0\le k\le m-1$.
Additionally, we shall say that $(f,g)\in C^\infty$ satisfy the
compatibility conditions to infinite order if this condition holds
for all $m$.

We can now state our main result:

\begin{theorem}\label{theorem1.1}
Let $\mathcal{K}$ be a fixed compact obstacle with smooth boundary
that satisfies \eqref{1.4}.  Assume that $Q(du,d^2u)$ and
$\Box$ are as above and that $(f,g)\in C^\infty(\ext)$
satisfy the compatibility conditions to infinite order.  Then
there is a constant $\varepsilon_0>0$, and an integer $N>0$ so
that for all $\varepsilon<\varepsilon_0$, if
\begin{equation}\label{1.11}
\sum_{|\alpha|\le N}\|<x>^{|\alpha|}\partial^\alpha_x
f\|_2+\sum_{|\alpha|\le N-1}\|<x>^{1+|\alpha|}\partial_x^\alpha
g\|_2 \le \varepsilon
\end{equation}
then \eqref{1.5} has a unique solution $u\in
C^\infty([0,\infty)\times\ext)$.
\end{theorem}

This paper is organized as follows.  In the next section, we will
recall some energy estimates from \cite{MS}.  In \S 3, we will
gather the pointwise estimates that we will require.  In \S 4, we will
collect some Sobolev-type estimates.  Included are some bounds on the
null forms which are exterior domain analogs of those from 
\cite{Si3} and 
\cite{So2}.
Finally, in \S 5, we will
use these estimates to prove the global existence theorem via a
continuity argument.


\newsection{$L^2$ Estimates.}
In this section, we will recall some estimates of \cite{MS}
related to the energy inequality.  Unless stated otherwise, the
proofs of the following results can be found in \cite{MS}.
Specifically, we will be concerned with solutions $u\in C^\infty
(\Rplus\times\ext)$ of the Dirichlet-wave equation
\begin{equation}\label{2.1}
\begin{cases}
\Box_\gamma u=F\\
u|_{\partial\mathcal{K}}=0\\
u|_{t=0}=f,\quad \partial_t u|_{t=0}=g
\end{cases}
\end{equation}
where
$$(\Box_\gamma u)^I=(\partial_t^2-c_I^2\Delta)u^I
+\sum_{J=1}^D\sum_{j,k=0}^3
\gamma^{IJ,jk}(t,x)\partial_j\partial_k u^J,\quad 1\le I\le D.$$
We shall assume that the $\gamma^{IJ,jk}$ satisfy the symmetry
conditions
\begin{equation}\label{2.2}
\gamma^{IJ,jk}=\gamma^{JI,jk}=\gamma^{IJ,kj}
\end{equation}
as well as the size condition
\begin{equation}\label{2.3}
\sum_{I,J=1}^D \sum_{j,k=0}^3 \|\gamma^{IJ,jk}(t,x)\|_{\infty}\le
\delta
\end{equation}
for $\delta$ sufficiently small (depending on the wave speeds).
 The energy estimate will involve bounds for the gradient of the
perturbation terms
$$\|\gamma'(t,\cd)\|_\infty = \sum_{I,J=1}^D\sum_{j,k,l=0}^3
\|\partial_l \gamma^{IJ,jk}(t,\cd)\|_\infty,$$ and the energy form
associated with $\Box_\gamma$, $e_0(u)=\sum_{I=1}^D e_0^I(u)$,
where
\begin{multline}\label{2.4}
e_0^I(u)=(\partial_0 u^I)^2+\sum_{k=1}^3 c_I^2(\partial_k u^I)^2
\\+2\sum_{J=1}^D\sum_{k=0}^3
\gamma^{IJ,0k}\partial_0u^I\partial_ku^J -
\sum_{J=1}^D\sum_{j,k=0}^3
\gamma^{IJ,jk}\partial_ju^I\partial_ku^J.
\end{multline}

The most basic estimate will lead to a bound for
$$E_M(t)=E_M(u)(t)=\int \sum_{j=0}^M e_0(\partial^j_t
u)(t,x)\:dx.$$

\begin{lemma}\label{lemma2.1}
Fix $M=0,1,2,\dots$, and assume that the perturbation terms
$\gamma^{IJ,jk}$ are as above.  Suppose also that $u\in C^\infty$
solves \eqref{2.1} and for every $t$, $u(t,x)=0$ for large $x$.
Then there is an absolute constant $C$ so that
\begin{equation}\label{2.5}
\partial_t E^{1/2}_M(t)\le C\sum_{j=0}^M \|\Box_\gamma
\partial_t^ju(t,\cd)\|_2+C\|\gamma'(t,\cd)\|_\infty E^{1/2}_M(t).
\end{equation}
\end{lemma}

Before stating the next result, let us introduce some notation. If
$P=P(t,x,D_t,D_x)$ is a differential operator, we shall let
$$[P,\gamma^{kl}\partial_k\partial_l]u=\sum_{1\le I,J\le
D}\sum_{0\le k,l\le 3} |[P,\gamma^{IJ,kl}\partial_k\partial_l]
u^J|.$$

In order to generalize the above energy estimate to include the
more general vector fields $L, Z$, we will need to use a variant
of the scaling vector field $L$.  We fix a bump function $\eta\in
C^\infty(\R^3)$ with $\eta(x)=0$ for $x\in \mathcal{K}$ and
$\eta(x)=1$ for $|x|>1$.  Then, set $\tilde{L}=\eta(x)r\partial_r
+ t\partial_t$. Using this variant of the scaling vector field and
an elliptic regularity argument, one can establish

\begin{proposition}\label{proposition2.2}
Suppose that the constant in \eqref{2.3} is small.  Suppose
further that
\begin{equation}\label{2.6}
\|\gamma'(t,\cd)\|_\infty \le \delta/(1+t),
\end{equation}
and
\begin{multline}\label{2.7}
\sum_{\substack{j+\mu\le
N_0+\nu_0\\\mu\le\nu_0}}\left(\|\tilde{L}^\mu\partial^j_t\Box_\gamma
u(t,\cd)\|_2+\|[\tilde{L}^\mu\partial_t^j,\gamma^{kl}\partial_k\partial_l]u(t,\cd)\|_2\right)\\
\le \frac{\delta}{1+t}\sum_{\substack{j+\mu\le N_0+\nu_0\\
\mu\le\nu_0}}\|\tilde{L}^\mu\partial_t^j
u'(t,\cd)\|_2+H_{\nu_0,N_0}(t),
\end{multline}
where $N_0$ and $\nu_0$ are fixed.  Then
\begin{multline}\label{2.8}
\sum_{\substack{|\alpha|+\mu\le N_0+\nu_0\\
\mu\le\nu_0}}\|L^\mu\partial^\alpha u'(t,\cd)\|_2 \\ \le
C\sum_{\substack{|\alpha|+\mu\le N_0+\nu_0-1\\
\mu\le\nu_0}}\|L^\mu\partial^\alpha\Box u(t,\cd)\|_2 +
C(1+t)^{A\delta}\sum_{\substack{\mu+j\le N_0+\nu_0\\ \mu\le
\nu_0}}\left(\int e_0(\tilde{L}^\mu\partial_t^j
u)(0,x)\:dx\right)^{1/2}\\
+C(1+t)^{A\delta}\Bigl(\int_0^t \sum_{\substack{|\alpha|+\mu\le
N_0+\nu_0-1\\\mu\le\nu_0-1}}\|L^\mu\partial^\alpha \Box
u(s,\cd)\|_2\:ds +\int_0^t H_{\nu_0,N_0}(s)\:ds\Bigr)\\
+C(1+t)^{A\delta}\int_0^t \sum_{\substack{|\alpha|+\mu\le
N_0+\nu_0\\ \mu\le \nu_0-1}}\|L^\mu\partial^\alpha
u'(s,\cd)\|_{L^2(|x|<1)}\:ds,
\end{multline}
where the constants $C$ and $A$ are absolute constants.
\end{proposition}

In practice $H_{\nu_0,N_0}(t)$ will involve weighted $L^2_x$ norms
of $|L^\mu \partial^\alpha u'|^2$ with $\mu+|\alpha|$ much smaller
than $N_0+\nu_0$, and so the integral involving $H_{\nu_0,N_0}$
can be dealt with using an inductive argument and weighted
$L^2_tL^2_x$ estimates that will be presented at the end of this
section.

In proving our existence results for \eqref{1.5}, the key step
will be to obtain a priori $L^2$-estimates involving $L^\mu
Z^\alpha u'$.  Begin by setting
\begin{equation}\label{2.9}
Y_{N_0,\nu_0}(t)=\int \sum_{\substack{|\alpha|+\mu\le
N_0+\nu_0\\\mu\le\nu_0}}e_0(L^\mu Z^\alpha u)(t,x)\:dx.
\end{equation}
We, then, have the following proposition which shows how the
$L^\mu Z^\alpha u'$ estimates can be obtained from the ones
involving $L^\mu\partial^\alpha u'$.
\begin{proposition}\label{proposition2.3}
Suppose that the constant $\delta$ in \eqref{2.3} is small and
that \eqref{2.6} holds.  Then,
\begin{multline}\label{2.10}
\partial_t Y_{N_0,\nu_0}\le C Y^{1/2}_{N_0,\nu_0} \sum_{\substack{
|\alpha|+\mu\le N_0+\nu_0\\ \mu\le\nu_0}} \|\Box_\gamma L^\mu Z^\alpha
u(t,\cd)\|_2 + C \|\gamma'(t,\cd)\|_\infty Y_{N_0,\nu_0} \\
+C \sum_{\substack{|\alpha|+\mu\le N_0+\nu_0+1\\ \mu\le
\nu_0}} \|L^\mu \partial^\alpha u'(s,\cd)\|^2_{L^2(|x|<1)}.
\end{multline}
\end{proposition}

As in \cite{KSS2} and \cite{KSS3} we shall also require some
weighted $L^2_tL^2_x$ estimates.  They will be used, for example,
to control the local $L^2$ norms such as the last term in
\eqref{2.10}.  For convenience, for the remainder of this section,
allow $\Box=\partial_t^2-\Delta$ to denote the unit speed
d'Alembertian.  The transition from the following estimates to
those involving \eqref{1.6} is straightforward.  Also, allow
$$S_T=\{[0,T]\times\ext\}$$
to denote the time strip of height $T$ in $\Rplus\times\ext$.

We, then, have the following proposition.
\begin{proposition}\label{proposition2.4}
Fix $N_0$ and $\nu_0$.  Suppose that $\mathcal{K}$ satisfies the
local exponential energy decay \eqref{1.4}.  Suppose further that
$u\in C^\infty$ solves \eqref{2.1} and $u(t,x)=0$ for $t<0$. Then
there is a constant $C=C_{N_0,\nu_0,\mathcal{K}}$ so that if $u$
vanishes for large $x$ at every fixed $t$
\begin{multline}\label{2.11}
(\log(2+T))^{-1/2}\sum_{\substack{|\alpha|+\mu\le
N_0+\nu_0\\\mu\le\nu_0}} \|<x>^{-1/2}L^\mu \partial^\alpha
u'\|_{L^2(S_T)} \\
\le C\int_0^T \sum_{\substack{|\alpha|+\mu\le
N_0+\nu_0+1\\\mu\le\nu_0}} \|\Box L^\mu \partial^\alpha
u(s,\cd)\|_2\:ds +C \sum_{\substack{|\alpha|+\mu\le
N_0+\nu_0\\\mu\le\nu_0}} \|\Box L^\mu\partial^\alpha
u\|_{L^2(S_T)}
\end{multline}
and
\begin{multline}\label{2.12}
(\log(2+T))^{-1/2}\sum_{\substack{|\alpha|+\mu\le
N_0+\nu_0\\\mu\le\nu_0}}\|<x>^{-1/2}L^\mu Z^\alpha
u'\|_{L^2(S_T)}\\
\le C\int_0^T \sum_{\substack{|\alpha|+\mu\le
N_0+\nu_0+1\\\mu\le\nu_0}} \|\Box L^\mu Z^\alpha u(s,\cd)\|_2\:ds
+ C \sum_{\substack{|\alpha|+\mu\le N_0+\nu_0\\\mu\le\nu_0}}
\|\Box L^\mu Z^\alpha u\|_{L^2(S_T)}.
\end{multline}
\end{proposition}

To be able to handle the last term in \eqref{2.8}, we shall need
the following.
\begin{lemma}\label{lemma2.5}
Suppose that \eqref{1.4} holds, and suppose that $u\in C^\infty$
solves \eqref{2.1} and satisfies $u(t,x)=0$ for $t<0$.  Then, for
fixed $N_0$ and $\nu_0$ and $t>2$,
\begin{multline}\label{2.13}
\sum_{\substack{|\alpha|+\mu\le N_0+\nu_0\\\mu\le\nu_0}} \int_0^t
\|L^\mu \partial^\alpha u'(s,\cd)\|_{L^2(|x|<2)}\:ds \\
\le C \sum_{\substack{|\alpha|+\mu\le N_0+\nu_0+1\\\mu\le\nu_0}}
\int_0^t \left(\int_0^s \|L^\mu \partial^\alpha \Box
u(\tau,\cd)\|_{L^2(||x|-(s-\tau)|<10)}\:d\tau\right)\:ds.
\end{multline}
\end{lemma}


\newsection{Pointwise Estimates.}
Here we will estimate solutions of the scalar inhomogeneous wave
equation
\begin{equation}\label{3.1}
\begin{cases}
(\partial_t^2-\Delta)w(t,x)=F(t,x),\quad (t,x)\in\R_+\times\ext\\
w(t,\cd)|_{\partial\mathcal{K}}=0\\
w(t,x)=0, \quad t\le 0.
\end{cases}
\end{equation}
If we assume, as before, that $\mathcal{K}\subset\{x\in\R^3\, :\,
|x|<1\}$ and that $\mathcal{K}$ satisfies \eqref{1.4}, then we have
the following result whose proof can be found in Metcalfe-Sogge
\cite{MS}.
\begin{theorem}\label{theorem3.1}
Let $w$ be a solution to \eqref{3.1}, and suppose that the
local energy decay bounds \eqref{1.4} hold for $\mathcal{K}$.  Then,
\begin{multline}\label{3.2}
(1+t+|x|)|L^\nu Z^\alpha w(t,x)|\le C\int_0^t \int_\ext
\sum_{\substack{|\beta|+\mu \le |\alpha|+\nu+7\\\mu\le\nu +1}} |L^\mu
Z^\beta F(s,y)|\:\frac{dy\:ds}{|y|}\\
+C\int_0^t \sum_{\substack{|\beta|+\mu\le |\alpha|+\nu+4\\ \mu\le\nu
+1}}\|L^\mu \partial^\beta F(s,\cd)\|_{L^2(|y|<2)}\:ds.
\end{multline}
\end{theorem}
\noindent Here and throughout $\{|y|<2\}$ is understood to mean
$\{y\in \ext\, : \, |y|<2\}$.

Additionally, we can prove the following improved pointwise bound for
the gradient of the solution $w$.  In this modified result, we are able to
bring the gradient inside the main term (the first term) on the right side.

\begin{theorem}\label{theorem3.2}
Let  $w$ be a solution to \eqref{3.1}, and suppose that the local
energy decay bounds \eqref{1.4} hold for $\mathcal{K}$.  Suppose
further that
$F(t,x)=0$ when $|x|>10 t$.  Then, if $|x|<t/10$ and $t>1$,
\begin{equation}\begin{split}\label{3.3}
(1+t+|x|)|L^\nu &Z^\alpha w'(t,x)|\le C\sum_{\substack{\mu+|\beta|\le \nu+|\alpha|+3 \\ \mu\le \nu+1} }\int_{t/100}^t \int_\ext
|L^{\mu}Z^{\beta}F'(s,y)|\:\frac{dy\:ds}{|y|}\\
&+ C \sup_{0\le s\le t}(1+s)\sum_{\substack{|\beta|+\mu\le
|\alpha|+\nu+4 \\ \mu\le\nu}}\|L^\mu Z^\beta F(s,\cd)\|_\infty\\
&+ C\sup_{0\le s\le t} (1+s)\sum_{\substack{|\beta|+\mu\le
|\alpha|+\nu+7\\ \mu\le \nu}}\int_0^s
\int_{\substack{||y|-(s-\tau)|\le 10\\ |y|\le
(600+\tau)/2}}|L^\mu Z^\beta F(\tau,y)|\:\frac{dy\:d\tau}{|y|}\\
&+C\sup_{0\le s\le t}\sum_{\substack{|\beta|+\mu\le |\alpha|+\nu+8\\
\mu\le\nu +1}}\int_{s/100}^s \int_{|y|\ge (1+\tau)/10}|L^\mu Z^\beta
F(\tau,y)|\:\frac{dy\:d\tau}{|y|}.
\end{split}\end{equation}
\end{theorem}

The remainder of this section will be dedicated to the proof of
\eqref{3.3}.

The Minkowski space estimate we shall use says that if $w_0$ is a
solution to the inhomogeneous wave equation
\begin{equation}\label{3.4}
\begin{cases}
(\partial_t^2-\Delta)w_0(t,x)=G(t,x),\quad (t,x)\in \R_+\times\R^3\\
w_0(0,x)=\partial_t w_0(0,x)=0,
\end{cases}
\end{equation}
and if $G(s,y)=0$ when $|y|>10s$, then
\begin{equation}\label{3.5}
(1+t)|L^\nu Z^\alpha w_0(t,x)|\le C\int_{t/100}^t
\int_{\R^3}\sum_{\substack{|\beta|+\mu \le |\alpha|+\nu+3\\ \mu\le \nu+1}}
|L^\mu Z^\beta G(s,y)|\:\frac{dy\:ds}{|y|}, \quad |x|<t/10.
\end{equation}
By using sharp Huygens principle, this follows from inequality
(2.3) in \cite{KSS3} and the fact that $[\partial_t^2-\Delta,Z]=0$ and
$[\partial_t^2-\Delta,L]=2(\partial_t^2 - \Delta)$.

Recall that we are assuming that $\mathcal{K}\subset \{x\in\R^3\, :\,
|x|<1\}$.  With this in mind, the first step is to see that
\eqref{3.5} yields for $|x|<t/10$
\begin{multline}\label{3.6}
(1+t)|L^\nu Z^\alpha w'(t,x)|\le C \int_{t/100}^t
\int_{\ext}\sum_{\substack{|\beta|+\mu\le |\alpha|+\nu+3\\\mu\le\nu
+1}}|L^\mu Z^\beta F'(s,y)|\:\frac{dy\:ds}{|y|}
\\+C\sup_{|y|\le 2,0\le s\le t}(1+s)\sum_{\substack{|\beta|+\mu\le
|\alpha|+\nu+1\\ \mu\le \nu}} |L^\mu \partial^\beta w'(s,y)|.
\end{multline}
The proof is exactly like that of Lemma 4.2 in \cite{KSS3}.

As a result of \eqref{3.6}, we would be done if we could show that
\begin{equation}\label{3.7}
\sup_{0\le s\le t}(1+s)\sum_{\substack{|\beta|+\mu\le
|\alpha|+\nu+1\\\mu\le\nu}}\|L^\mu \partial^\beta w'(s,\cd)\|_{L^\infty(|x|<2)}
\end{equation}
is controlled by the right side of \eqref{3.3}.

To prove this, we shall need the following
\begin{lemma}\label{lemma3.3}
Suppose that $w$ is as above.  Suppose further that $(\partial_t^2-\Delta)
w(s,y)=F(s,y)=0$ if $|y|>10$.  Then,
\begin{equation}\label{3.8}
(1+t)\sup_{|x|<2}|L^\nu \partial^\alpha w'(t,x)|\le C\sup_{0\le s\le
t}\sum_{\substack{|\beta|+\mu\le \nu+|\alpha|+3\\ \mu\le \nu}}(1+s)\|L^\mu
\partial^\beta F(s,\cd)\|_2.
\end{equation}
\end{lemma}

\noindent{\bf Proof of Lemma \ref{lemma3.3}:} By Sobolev estimates,
the left side of \eqref{3.8} is dominated by
$$(1+t)\sum_{\substack{|\beta|+\mu\le \nu+|\alpha|+2\\ \mu\le \nu}}\|L^\mu
\partial^\beta w'(t,\cd)\|_{L^2(|x|<3)}.$$
By exponential energy decay and elliptic regularity (see Lemma 2.8 in \cite{MS}), there must be a
constant $c>0$ so that this is controlled by
\begin{multline*}
(1+t)\sum_{\substack{|\beta|+\mu\le
\nu+|\alpha|+1\\\mu\le\nu}}\|L^\mu\partial^\beta F(t,\cd)\|_2+(1+t)\int_0^t
e^{-c(t-s)} \sum_{\substack{|\beta|+\mu\le \nu+|\alpha|+3\\\mu\le\nu}}\|L^\mu
\partial^\beta F(s,\cd)\|_2\:ds\\
\le C\sup_{0\le s\le t}(1+s)\sum_{\substack{|\beta|+\mu\le \nu+
|\alpha|+3\\\mu\le\nu}}\|L^\mu \partial^\beta F(s,\cd)\|_2,
\end{multline*}
as desired.\qed

We also need an estimate for solutions whose forcing terms vanish near
the obstacle.  Suppose that $w$ is as above, but now assume that
$(\partial_t^2-\Delta)w(s,y)=F(s,y)=0$ when $|y|<5$.  Then, write
$$w=w_0+w_r$$
where $w_0$ solves the boundaryless wave equation
$(\partial_t^2-\Delta)w_0=F$ with zero initial data.  Fix $\eta\in
C_0^\infty(\R^3)$ satisfying $\eta(x)=1$ for $|x|<2$ and $\eta(x)=0$
for $|x|\ge 3$.  If we set $\tilde{w}=\eta w_0+w_r$, then since $\eta
F=0$, $\tilde{w}$ solves the Dirichlet-wave equation
$$(\partial_t^2-\Delta)\tilde{w}=G=-2\nabla_x \eta\cdot \nabla_x w_0 -
(\Delta \eta) w_0$$
with zero initial data.  This forcing term vanishes unless $2\le
|x|\le 3$.  In this case, by Lemma \ref{lemma3.3},
\begin{align*}
(1+t) \sum_{\substack{|\beta|+\mu \le |\alpha|+\nu+1 \\ \mu\le \nu}} \sup_{|x|<2}&|L^\mu \partial^\beta w'(t,x)| 
=(1+t)\sum_{\substack{|\beta|+\mu \le |\alpha|+\nu+1 \\ \mu\le \nu}} \sup_{|x|<2}|L^\mu \partial^\beta \tilde{w}'(t,x)|\\
&\le C \sup_{0\le s\le t}\sum_{\substack{|\beta|+\mu\le |\alpha|+\nu+4\\
\mu\le \nu}}(1+s)\|L^\mu \partial^\beta G(s,\cd)\|_2\\
&\le C \sup_{0\le s\le t}\sum_{\substack{|\beta|+\mu\le |\alpha|+\nu+5\\
\mu\le\nu}} (1+s)\|L^\mu \partial^\beta w_0(s,\cd)\|_{L^\infty(2<|x|<3)}.
\end{align*}

Thus, if we could show the following lemma, we would have that the
\eqref{3.7} is bounded by the right side of \eqref{3.3}, which
would complete the proof of Theorem \ref{theorem3.2}.

\begin{lemma}\label{lemma3.4}
Suppose that $v$ is a solution to the free wave equation
$(\partial_t^2-\Delta)v=G$ and that $v$ has vanishing Cauchy data.
Suppose further that $G(t,x)=0$ when $|x|>10t$.  Then,
\begin{multline}\label{3.9}
\sup_{2<|x|<3} |v(t,x)|\le C \sum_{|\alpha|\le 2} \int_0^t
\int_{\substack{||y|-(t-s)| <10\\ |y|\le (600+s)/2}}
|\Omega^\alpha G(s,y)|\:\frac{dy\:ds}{|y|} \\
+\frac{C}{1+t}\sum_{\substack{|\alpha|+\mu \le 3 \\ \mu\le 1}}\int_{t/100}^t \int_{|y|\ge
(1+s)/10} |L^\mu Z^\alpha G(s,y)|\:\frac{dy\:ds}{|y|}.
\end{multline}
\end{lemma}

\noindent {\bf Proof of Lemma \ref{lemma3.4}:}  This lemma is a
consequence of the estimate
\begin{equation}\label{3.10}
|x||v(t,x)|\le \frac{1}{2}\int_0^t \int_{|r-(t-s)|}^{r+t-s}
 \sup_{|\theta|=1} |G(s,\rho\theta)|\rho\:d\rho\:ds,
\end{equation}
where $r=|x|$.  See, e.g., (2.4) in \cite{KSS3}.

We begin by choosing a cutoff
function $\rho(x)$ satisfying $\rho(x)=1$ for $|x|<1/10$ and
$\rho(x)=0$ for $|x|>1/2$.  Set $G_1(t,x)=\rho(x/(1+t))G(t,x)$ and
$G_2(t,x)=(1-\rho(x/(1+t)))G(t,x)$, and for $j=1,2$, let $v_j$ solve the
inhomogeneous wave equation $(\partial_t^2-\Delta)v_j=G_j$ with zero
initial data.  Then, $v=v_1+v_2$.  Using \eqref{3.10} and the
Sobolev estimate on the sphere, we have for $2<|x|<3$
\begin{align*}
|v_1(t,x)|&\le C\sum_{|\alpha|\le 2}\int_0^t
 \int_{||y|-(t-s)|<10} |\Omega^\alpha
 G_1(s,y)|\:\frac{dy\:ds}{|y|}\\
&\le C\sum_{|\alpha|\le 2} \int_0^t
 \int_{\substack{||y|-(t-s)|<10\\ |y|\le (1+s)/2}}
 |\Omega^\alpha G(s,y)|\:\frac{dy\:ds}{|y|}.
\end{align*}

By \eqref{3.5}, if $2<|x|<3$ and $t>10|x|$,
\begin{align*}
|v_2(t,x)|&\le \frac{C}{1+t}\sum_{\substack{|\alpha|+\mu\le 3 \\ \mu\le
 1}}\int_{t/100}^t \int |L^\mu Z^\alpha G_2(s,y)|\:\frac{dy\:ds}{|y|}\\
&\le \frac{C}{1+t}\sum_{\substack{|\alpha|+\mu\le 3 \\ \mu\le 1}}\int_{t/100}^t
 \int_{|y|\ge (1+s)/10}|L^\mu Z^\alpha G(s,y)|\:\frac{dy\:ds}{|y|}.
\end{align*}

Since $v=v_1+v_2$, these two estimates yield \eqref{3.9} when
$t>10|x|$.  The proof of the estimate for $v_1$ shows that for $0<t<10|x|$
the left side of \eqref{3.9} is dominated by the first term in the right.\qed


\newsection{Estimates Related to the Null Condition and Sobolev-type Estimates.}
The first result of this section concerns bounds for the null forms.
They must involve the weight $\langle c_Jt-r\rangle$ since we are not
using the generators of Lorentz rotations.  The estimates will involve
the admissible homogeneous vector fields that we are using
$\{\Gamma\}=\{Z,L\}$.  The proof of
these estimates can be found in Sideris-Tu \cite{Si3} and
Sogge \cite{So2}.

\begin{lemma}\label{lemma4.1}
Suppose that the quasilinear null condition \eqref{1.9} holds.
Then,
\begin{multline}\label{4.1}
\Bigl|\sum_{0\le j,k,l\le 3} B^{JJ,jk}_{J,l}\partial_l
u\partial_j\partial_k v\Bigr|\\
\le C \langle r\rangle^{-1}(|\Gamma u||\partial^2 v|+|\partial
u||\partial\Gamma v|) + C \frac{\langle c_Jt-r\rangle}{\langle
t+r\rangle}|\partial u||\partial^2 v|.
\end{multline}
Also, if the semilinear null condition \eqref{1.10}
holds
\begin{equation}\label{4.2}
\Bigl|\sum_{0\le j,k\le 3} A^{J,jk}_{JJ}\partial_j u \partial_k
v\Bigr| \le C \langle r\rangle^{-1}(|\Gamma u||\partial
v|+|\partial u||\Gamma v|)+ C\frac{\langle c_Jt-r\rangle}{\langle
t+r\rangle} |\partial u||\partial v|.
\end{equation}
\end{lemma}

We shall also need the following Sobolev-type estimate.  The first is
an exterior domain analog of results of Klainerman-Sideris \cite{KS}.

\begin{lemma}\label{lemma4.2}  Suppose that $u(t,x)\in
C^\infty_0(\mathbb{R}\times \mathbb{R}^3 \backslash \mathcal{K})$
vanishes for $x\in \partial\mathcal{K}$.  Then if $|\alpha|=M$ and
$\nu$ are fixed
\begin{multline}\label{4.3}
\|\langle t-r\rangle L^\nu Z^\alpha \partial^2 u(t,\cd)\|_2\le
C\sum_{\substack{|\beta|+\mu \le M+\nu + 1 \\ \mu \le
\nu+1}}\|L^\mu Z^\beta u'(t,\cd)\|_2 \\
+C\sum_{\substack{|\beta|+\mu \le M+\nu \\ \mu\le \nu}} \|\langle
t+r\rangle L^\mu Z^\beta (\partial_t^2-\Delta) u(t,\cd)\|_2 +
C(1+t)\sum_{\mu\le \nu} \|L^\mu u'(t,\cd)\|_{L^2(|x|<2)}.
\end{multline}
\end{lemma}

\noindent{\bf Proof of Lemma \ref{lemma4.2}:}  The first step is to show that
\begin{align}\label{4.4}
\|\langle t-r\rangle L^\nu Z^\alpha \partial^2 u(t,\cd)\|_2 &\le
C\sum_{\substack{|\beta|+\mu \le M+\nu +1 \\ \mu\le \nu+1}}\|L^\mu
Z^\beta u'(t,\cd)\|_2
\\
&+C\sum_{\substack{|\beta|+\mu\le M+\nu \\
\mu\le \nu}}\|\langle t+r\rangle L^\mu Z^\beta
(\partial^2_t-\Delta)u(t,\cd)\|_2 \notag
\\
&+C(1+t)\sum_{\substack{|\beta|+\mu\le M+\nu+2 \\ \mu\le
\nu}}\|L^\mu \partial^\beta u(t,\cd)\|_{L^2(|x|<3/2)} \notag
\end{align}

If one replaces the left side by the analogous expression with the
norm taken over $|x|<3/2$ then this term is dominated by the last
term in \eqref{4.4} due to the fact that the coefficients of $Z$
are bounded when $|x|<3/2$.

To handle the part where $|x|>3/2$ we shall use the following
Minkowski space estimate
\begin{multline}\label{4.5}
\|\langle t-r\rangle L^\nu Z^\alpha \partial^2
h(t,\cd)\|_{L^2(\mathbb{R}^3)} \le C\sum_{\substack{|\beta|+\mu
\le |\alpha|+\nu+1 \\ \mu\le \nu+1}}\|L^\mu Z^\beta
h'(t,\cd)\|_{L^2(\mathbb{R}^3)} \\
+C\sum_{\substack{|\beta|+\mu \le |\alpha|+\nu \\ \mu\le
\nu}}\|\langle t+r\rangle L^\mu Z^\beta
(\partial^2_t-\Delta)h(t,\cd)\|_{L^2(\mathbb{R}^3)},
\end{multline}
which is valid for $h\in C^\infty_0(\mathbb{R}\times
\mathbb{R}^3)$.  This estimate follows from (2.10) and Lemma 3.1
of Klainerman and Sideris \cite{KS} if one uses the fact that
$[(\partial^2_t-\Delta), Z]=0$ and
$[(\partial^2_t-\Delta),L]=2(\partial_t-\Delta)$.

To use this, choose $\eta\in C^\infty_0(\mathbb{R}^3)$ so that
$\eta(x)=0$ for $|x|<1$ and $\eta(x)=1$ for $|x|>3/2$.  Then if we
let $h(t,x)=\eta(x)u(t,x)$, we have
$$(\partial^2_t-\Delta)h=\eta(x)(\partial^2_t-\Delta)u-2\nabla_x\eta\cdot\nabla_x
u-(\Delta \eta)u.$$ Therefore since the last two terms are
supported in $|x|<3/2$, \eqref{4.5} yields
\begin{align*}
\|\langle t-r\rangle L^\nu Z^\alpha \partial^2
u(t,\cd)\|_{L^2(|x|>3/2)} &\le \|\langle t-r\rangle L^\nu Z^\alpha
\partial^2 h(t,\cd)\|_{L^2(\mathbb{R}^3)}
\\
&\le C\sum_{\substack{|\beta|+\mu \le M+\nu+1 \\ \mu \le
\nu+1}}\|L^\mu Z^\beta u'(t,\cd)\|_2
\\
&+ C\sum_{\substack{|\beta|+\mu \le M+\nu \\ \mu \le \nu}}\|\langle
t+r\rangle L^\mu Z^\beta (\partial^2_t-\Delta)u(t,\cd)\|_2
\\
&+C(1+t)\sum_{\substack{|\beta|+\mu \le M+\nu + 1 \\ \mu\le
\nu}}\|L^\mu \partial^\beta u(t,\cd)\|_{L^2(|x|<3/2)},
\end{align*}
which completes the proof of \eqref{4.4}.

In view of this inequality, to finish the proof of the lemma we
need to show that if $1<R<2$ then
\begin{equation}\label{4.6}
(1+t)\sum_{\substack{|\beta|+\mu \le M+\nu+2 \\ \mu\le
\nu}}\|L^\mu \partial^\beta u(t,\cd)\|_{L^2(|x|<R)}
\end{equation}
is controlled by the right hand side of \eqref{4.3}.  However, if
$1<R<R_0<2$, by elliptic regularity this term is dominated by
\begin{multline*}
(1+t)\Bigl(\sum_{\substack{|\beta|+\mu \le M+\nu \\ \mu\le
\nu}}\|L^\mu
\partial^\beta \Delta u(t,\cd)\|_{L^2(|x|<R_0)} +
\sum_{\substack{|\beta|+\mu \le M+\nu  \\ \mu\le \nu}}\|L^\mu
\partial^\beta \partial_t u'(t,\cd)\|_{L^2(|x|<R)}
\\
+\sum_{\substack{|\beta|+\mu\le (M-1)+\nu +2 \\ \mu \le
\nu}}\|L^\mu \partial^\beta u(t,\cd)\|_{L^2(|x|<R_0)}\Bigr).
\end{multline*}
Since Lemma 2.3 from \cite{KS} yields
\begin{multline*}
\langle t-r\rangle\Bigl(\sum_{\substack{|\beta|+\mu \le M+\nu \\
\mu\le \nu}}|L^\mu
\partial^\beta \Delta u|+\sum_{\substack{|\beta|+\mu \le M+\nu
\\ \mu\le \nu }}|L^\mu \partial^\beta \partial_t u'|\Bigr)
\\
\le C\sum_{\substack{|\beta|+\mu \le M+\nu +1 \\ \mu\le
\nu+1}}|L^\mu Z^\beta u'| + C\langle t+r\rangle
\sum_{\substack{|\beta|+\mu \le M+\nu \\ \mu\le \nu}}|L^\mu
Z^\beta (\partial_t^2-\Delta)u|,
\end{multline*}
the preceding inequality and an induction argument imply that
\eqref{4.6} is dominated by the right hand side of \eqref{4.3}
plus
$$(1+t)\sum_{|\beta|\le 1, \mu\le \nu}\|L^\mu \partial^\beta u(t,\cd)\|_{L^2(|x|<\tilde{R})},$$
where $\tilde{R}$ can be taken to satisfy $3/2<\tilde{R}<2$.  Since
$\partial^j_tu$ vanishes on $\partial \mathcal{K}$ 
one can use a similar induction argument to see that this is also
dominated by the right hand side of \eqref{4.3} plus
$$(1+t)\sum_{\mu\le \nu}\|L^\mu \nabla_x
u(t,\cd)\|_{L^2(|x|<R_1)},$$ where $\tilde{R}<R_1<2$, which finishes the
proof. \qed

The next lemma is an exterior domain analog of an estimate of
Hidano-Yokoyama \cite{HY}.
\begin{lemma}\label{lemma4.3}
Suppose that $u(t,x)\in C^\infty_0(\R\times\ext)$ vanishes for
$x\in\partial\mathcal{K}$.  Then
\begin{multline}\label{4.7}
r^{1/2}\langle t-r\rangle |\partial L^\nu Z^\alpha  u(t,x)|\le
C\sum_{\substack{|\beta|+\mu\le |\alpha|+\nu+2\\ \mu\le
\nu+1}}\|L^\mu Z^\beta u'(t,\cd)\|_2 \\+
C\sum_{\substack{|\beta|+\mu\le |\alpha|+\nu+1\\\mu\le\nu}}
\|\langle t+r\rangle L^\mu Z^\beta (\partial_t^2-\Delta)
u(t,\cd)\|_2 +C(1+t)\sum_{\mu\le \nu} \|L^\mu
u'(t,\cd)\|_{L^\infty(|x|<2)}.
\end{multline}
\end{lemma}

\noindent{\bf Proof of Lemma \ref{lemma4.3}:}  Inequality (4.2) of
Hidano-Yokoyama \cite{HY} implies that in Min\-kow\-ski space
\begin{multline}\label{4.8}
r^{1/2}\langle t-r\rangle |\partial L^\nu Z^\alpha h(t,x)| 
\\\le C
\sum_{\substack{|\beta|+\mu\le |\alpha|+\nu+1 \\ \mu\le \nu}}
\Bigl(\|L^\mu Z^\beta h'(t,\cd)\|_2 + \|\langle t-r\rangle L^\mu
Z^\beta
\partial^2 h(t,\cd)\|_2\Bigr).
\end{multline}
If we choose $\eta\in C_0^\infty(\R^3)$ so that $\eta(x)=0$ for
$|x|<1$ and $\eta(x)=1$ for $|x|>5/4$, and let $h(t,x)=\eta(x)u(t,x)$,
then we conclude that when $|x|>5/4$,
\begin{align*}
r^{1/2}&\langle t-r\rangle |\partial L^\nu Z^\alpha  u(t,x)|\\
&\le C
\sum_{\substack{|\beta|+\mu\le |\alpha|+\nu+1 \\ \mu\le \nu}}
\Bigl(\|L^\mu Z^\beta u'(t,\cd)\|_2+\|\langle t-r\rangle  L^\mu
Z^\beta\partial^2
u(t,\cd)\|_2\Bigr) \\
&\quad\quad\quad\quad\quad\quad\quad\quad\quad\quad\quad\quad\quad+
C(1+t)\sum_{\substack{|\beta|+\mu \le |\alpha|+\nu+3 \\ \mu\le
\nu}} \|L^\mu \partial^\beta u(t,\cd)\|_{L^2(|x|<3/2)}.
\end{align*}
By the Sobolev inequality, over $|x|<5/4$ the left side of
\eqref{4.7} is bounded by a similar inequality involving only 
the last term on the right.

%

If we use (4.3), we see that the second term in the right is
dominated by the right side of \eqref{4.7}.  If we repeat
the last part of the proof of Lemma \ref{lemma4.2}, we conclude that the same
is true for the last term in the preceding inequality.  \qed

Finally, we will need the following now standard consequence of
the Sobolev lemma (see \cite{knull}).

\begin{lemma}\label{sobolev}
Suppose that $h\in C^\infty(\R^3)$.  Then, for $R>2$,
$$\|h\|_{L^\infty(R<|x|<R+1)}\le C R^{-1}\sum_{|\alpha|+|\beta|\le
2}\|\Omega^{\alpha}\partial^\beta_x h\|_{L^2(R-1<|x|<R+2)}.$$
\end{lemma}


\newsection{Global Existence and the Continuity Argument.}
In this section, we will prove the main result, Theorem
\ref{theorem1.1}.  We shall take $N=101$ in its smallness
hypothesis \eqref{1.11}, but this is certainly not optimal.

To prove our global existence theorem, we shall need a standard
local existence theorem.

\begin{theorem}\label{theorem5.1}
Suppose that $f$ and $g$ are as in Theorem \ref{theorem1.1} with
$N\ge 6$ in \eqref{1.11}.  Then there is a $T>0$ so that the
initial value problem \eqref{1.5} with this initial data has a
$C^2$ solution satisfying
$$u\in L^\infty([0,T];H^N(\ext))\cap C^{0,1}([0,T];
H^{N-1}(\ext)).$$ The supremum of such $T$ is equal to the
supremum of all $T$ such that the initial value problem has a
$C^2$ solution with $\partial^\alpha u$ bounded for $|\alpha|\le
2$.  Also, one can take $T\ge 2$ if $\|f\|_{H^N}+\|g\|_{H^{N-1}}$
is sufficiently small.
\end{theorem}

This essentially follows from the local existence results Theorem
9.4 and Lemma 9.6 in \cite{KSS}.  The latter were only stated for
diagonal single-speed systems; however, since the proof relied
only on energy estimates, it extends to the multi-speed
non-diagonal case if the symmetry assumptions \eqref{1.8} are
satisfied.

Next, as in \cite{KSS3}, in order to avoid dealing with
compatibility conditions for the Cauchy data, it is convenient to
reduce the Cauchy problem \eqref{1.5} to an equivalent equation
with a nonlinear driving force but vanishing Cauchy data.  We will
then set up a continuity argument that utilizes the results of the
previous three sections to show global existence and prove Theorem
\ref{theorem1.1}.

Recall that our smallness condition on the data is
\begin{equation}\label{5.1}
\sum_{|\alpha|\le 101} \|<x>^{|\alpha|}\partial_x^\alpha
f\|_{L^2(\ext)}+\sum_{|\alpha|\le
100}\|<x>^{1+|\alpha|}\partial_x^\alpha g\|_{L^2(\ext)}\le
\varepsilon.
\end{equation}
To make the reduction to an equation with vanishing initial data,
we will begin by noting that if the data satisfies \eqref{5.1} for
$\varepsilon$ sufficiently small, then we can find a solution $u$
to \eqref{1.5} on a set of the form $0<ct<|x|$ where $c=5\max_I
c_I$, and this solution satisfies
\begin{equation}\label{5.2}
\sup_{0<t<\infty}\sum_{|\alpha|\le
101}\|<x>^{|\alpha|}\partial^{\alpha}u(t,\cd)\|_{L^2(\ext :
|x|>ct)}\le C_0\varepsilon,
\end{equation}
for an absolute constant $C_0$.

To prove this, we shall repeat the argument of Keel-Smith-Sogge
\cite{KSS3}.  By scaling in $t$, we may assume without loss that
$\max_I c_I=1/2$.  Theorem \ref{theorem5.1} yields a solution $u$
to \eqref{1.5} on the set $0<t<2$ which satisfies \eqref{5.2}.  We
wish to show that this solution extends to the set $0<ct<|x|$.  To
do so, let $R\ge 4$ and consider data $(f_R,g_R)$ supported in the
set $R/4<|x|<4R$ which agrees with the data $(f,g)$ on the set
$R/2<|x|<2R$.  Let $u_R(t,x)$ satisfy the free wave equation
$$\Box u_R=Q(du_R, R^{-1}d^2u_R)$$
with Cauchy data $(f_R(R\cd), Rg_R(R\cd))$.  The solution $u_R$
then exists for $0<t<1$ by standard local existence theory (see,
e.g., \cite{H} and \cite{S}) and satisfies
\begin{multline*}
\sup_{0<t<1} \|u_R(t,\cd)\|_{H^{101}(\R^3)}\le
C(\|f_R(R\cd)\|_{H^{101}(\R^3)}+R\|g_R(R\cd)\|_{H^{100}(\R^3)})\\
\le C R^{-3/2} \Bigl(\sum_{|\alpha|\le 101} \|(R\partial_x)^\alpha
f_R\|_{L^2(\R^3)}+R\sum_{|\alpha|\le 100}\|(R\partial_x)^\alpha
g_R\|_{L^2(\R^3)}\Bigr).
\end{multline*}
The smallness condition on $|u'_R|$ implies that the wave speeds
for the quasilinear equation are bounded above by 1.  A domain of
dependence argument shows that the solutions
$u_R(R^{-1}t,R^{-1}x)$ restricted to $||x|-R|<\frac{R}{2}-t$ agree
on their overlaps, and also with the local solution, yielding a
solution to \eqref{1.5} on the set $\{\ext : 2t<|x|\}$.  An
argument using a partition of unity now yields \eqref{5.2}.

We are now ready to set up the continuity argument.  We will use
the local solution $u$ to allow us to restrict to the case where
the Cauchy data vanish.  Fix a cutoff function $\chi\in
C^\infty(\R)$ satisfying $\chi(s)=1$ if $s\le \frac{1}{2c}$ and
$\chi(s)=0$ if $s>\frac{1}{c}$.  Set
$$u_0(t,x)=\eta(t,x)u(t,x), \quad \eta(t,x)=\chi(|x|^{-1}t).$$
Assuming as we may that $0\in\mathcal{K}$, we have that $|x|$ is
bounded below on the complement of $\mathcal{K}$ and the function
$\eta(t,x)$ is smooth and homogeneous of degree 0 in $(t,x)$.
Additionally,
$$\Box u_0 = \eta Q(du, d^2u)+[\Box, \eta] u.$$
Thus, $u$ solves $\Box u=Q(du,d^2u)$ for $0<t<T$ if and only if
$w=u-u_0$ solves
\begin{equation}\label{5.3}
\begin{cases}
\Box w = (1-\eta)Q(du,d^2u)-[\Box, \eta] u\\
w|_{\partial\mathcal{K}}=0\\
w(t,x)=0,\quad t\le 0
\end{cases}
\end{equation}
for $0<t<T$.

A key step in proving that \eqref{5.3} admits a global solution is
to prove uniform energy and dispersive estimates for $w$ on the
interval of existence.  First note that since $u_0=\eta u$, by
\eqref{5.2} and Lemma \ref{sobolev}, there is an absolute constant
$C_1$ so that
\begin{multline}\label{5.4}
(1+t+|x|)\sum_{\mu+|\alpha|\le 99}|L^\mu Z^\alpha u_0(t,x)|\\
+\sum_{\mu+|\alpha|+|\beta|\le 101} \|<t+r>^{|\beta|}L^\mu
Z^\alpha \partial^\beta u_0(t,\cd)\|_2\le C_1 \varepsilon.
\end{multline}
Furthermore, if we let $v$ be the solution of the linear equation
\begin{equation}\label{5.5}
\begin{cases}
\Box v=-[\Box,\eta]u\\
v|_{\partial\mathcal{K}}=0\\
v(t,x)=0,\quad t\le 0,
\end{cases}
\end{equation}
then we will show that \eqref{5.2} and Theorem \ref{theorem3.1}
imply that there is an absolute constant $C_2$ so that
\begin{equation}\label{5.6}
(1+t+|x|)\sum_{\mu+|\alpha|\le 92}|L^\mu Z^\alpha v(t,x)| +
\sum_{\substack{\mu+|\alpha|\le 90}} \|L^\mu Z^\alpha
v'(t,\cd)\|_2 \le C_2 \varepsilon.
\end{equation}

Indeed, by \eqref{3.2}, the first term on the left side of
\eqref{5.6} is bounded by
\begin{multline*}
\int_0^t \int_{|x|>cs} \sum_{\mu+|\alpha|\le 99}|L^\mu Z^\alpha
([\Box, \eta]u)(s,x)|\:\frac{dx\:ds}{|x|} \\
+\int_0^t \sum_{\mu+|\alpha|\le 96} \|L^\mu \partial^\beta ([\Box,
\eta]u)(s,\cd)\|_{L^2(\ext : |x|<2)}\:ds
\end{multline*}
which by the Schwarz inequality is bounded by
$$\sum_{\mu+|\alpha|\le 99}\sum_{j=0}^\infty \sup_{0<cs<2^{j+1}}
\|<x>^{3/2} L^\mu Z^{\alpha} [\Box, \eta] u(s,\cd)\|_{L^2(\ext :
2^j<|x|<2^{j+1})}.$$ Since this is dominated by
$$\sup_{0<t<\infty}\sum_{\mu+|\alpha|\le 99}\|<x>^2 L^\mu Z^\alpha
[\Box, \eta] u(t,\cd)\|_2,$$ one gets that the first term on the
left side of \eqref{5.6} is $O(\varepsilon)$ from \eqref{5.2} and
the homogeneity of $\eta$.

For the second term on the left side of \eqref{5.6}, if we argue
as in the proof of \eqref{2.5} (except now for the linear wave
equation), we see that
\begin{multline*}
\partial_t \sum_{\substack{\mu+|\alpha|\le 90}} \|L^\mu
Z^\alpha v'(t,\cd)\|^2_2 \\\le C
\Bigl(\sum_{\substack{\mu+|\alpha|\le 90}} \|L^\mu Z^\alpha
 v'(t,\cd)\|_2\Bigr)
\Bigl(\sum_{\substack{\mu+|\alpha|\le 90}} \|L^\mu Z^\alpha \Box
v(t,\cd)\|_2\Bigr) \\
+ C\sum_{\substack{\mu+|\alpha|\le 90}}\Bigl|\int_{\partial\mathcal{K}}
 \partial_0 L^\mu Z^\alpha
v(t,\cd)\nabla L^\mu Z^\alpha v(t,\cd)\cdot n\:d\sigma\Bigr|,
\end{multline*}
where $n$ is the outward normal at a given point on
$\partial\mathcal{K}$. 
Since $\mathcal{K}\subset\{|x|<1\}$, it
follows that
\begin{multline*}
\partial_t \sum_{\substack{\mu+|\alpha|\le 90}} \|L^\mu
Z^\alpha v'(t,\cd)\|^2_2 \\\le C
\Bigl(\sum_{\substack{\mu+|\alpha|\le 90}} \|L^\mu Z^\alpha
 v'(t,\cd)\|_2\Bigr)
\Bigl(\sum_{\substack{\mu+|\alpha|\le 90}} \|L^\mu Z^\alpha \Box
v(t,\cd)\|_2\Bigr) \\+C \int_{\{x\in\ext : |x|<1\}}
\sum_{\substack{\mu+|\alpha|\le 91}} |L^\mu Z^\alpha
v'(t,\cd)|^2\:dx.
\end{multline*}
Thus, since $\Box v(s,y)=-[\Box,\eta]u(s,y)$, it follows that
\begin{multline}\label{5.7}
\sum_{\mu+|\alpha|\le 90}\|L^\mu Z^\alpha v'(t,\cd)\|^2_2
\le C\Bigl(\int_0^t \sum_{\mu+|\alpha|\le 90}\|L^\mu Z^\alpha
(-[\Box,\eta]u)(s,y)\|_2\:ds \Bigr)^2
\\+ C\int_0^t \sum_{\mu+|\alpha|\le
91}\|L^\mu Z^\alpha v'(s,\cd)\|^2_{L^2(|x|<1)}\:ds.
\end{multline}
The first term on the right is $O(\varepsilon)$ by \eqref{5.2}.  Using
the bound for the first term in the left of \eqref{5.6}, it follows
that the second term on the right of \eqref{5.7} is also
$O(\varepsilon)$ as desired.

Using this, we are now ready to set up the continuity argument. If
$\varepsilon>0$ is as above, we shall assume that we have a
solution of our equation \eqref{1.5} for $0\le t\le T$ such that
we have the following estimates
\begin{align}\label{5.8}
\sum_{\substack{|\alpha|+\nu\le 52\\\nu\le 2}} \|L^\nu Z^\alpha
w'(t,\cdot)\|_2 &\le A_0\varepsilon\\
(1+t+r)\sum_{|\alpha|\le 40}|Z^\alpha w'(t,x)|& \le A_1\varepsilon\label{5.9}\\
(1+t+r)\sum_{\substack{|\alpha|+\nu\le 55\\ \nu\le 3}} |L^\nu
Z^\alpha (w-v)(t,x)|&\le B_1 \varepsilon^2 (1+t)^{1/10}\log(2+t) \label{5.10}\\
\sum_{|\alpha|\le 100} \|\partial^\alpha u'(t,\cd)\|_2 &\le B_2 \varepsilon
(1+t)^{1/40}\label{5.11}\\
\sum_{\substack{|\alpha|+\nu\le 65\\ \nu\le 4}}
\|L^\nu Z^\alpha u'(t,\cd)\|_2&\le B_3 \varepsilon (1+t)^{1/20} \label{5.12}\\
\sum_{\substack{|\alpha|+\nu\le 63\\ \nu\le 4}}\|\langle
x\rangle^{-1/2}L^\nu Z^\alpha u'\|_{L^2(S_t)}&\le B_4\varepsilon
(1+t)^{1/20}(\log (2+t))^{1/2}. \label{5.13}
\end{align}
Here, as before, the $L^2$ norms are taken over $\ext$ and the
weighted $L^2_tL^2_x$ norms are taken over $S_t=[0,t]\times\ext$.
In the main estimates \eqref{5.8} and \eqref{5.9}, we can take
$A_0=A_1=4C_2$, where $C_2$ is the constant occurring in the
bounds \eqref{5.6} for $v$.

Clearly if $\varepsilon$ is small then all of these estimates are
valid, if $T=2$, by Theorem \ref{theorem5.1}. Keeping this in
mind, we shall then prove that, for $\varepsilon>0$ smaller than
some number depending on $B_1, \dots, B_4$,
\begin{enumerate}
\item[\bf i.)] \eqref{5.8} is valid with $A_0$ replaced by
$A_0/2$. \item[\bf ii.)] Under the assumption of {\bf{(i.)}}, that
\eqref{5.9} is valid with $A_1$ replaced by $A_1/2$.  \item[\bf
iii.)] \eqref{5.10}-\eqref{5.13} are  consequences of \eqref{5.8}
and \eqref{5.9} for suitable constants $B_i$.
\end{enumerate}
By the local existence theorem, it will follow that a solutions
exists for all $t>0$ if $\varepsilon$ is small enough.

Before we begin the proof of (i.), we will set up some
preliminary results under the assumption of \eqref{5.8}-\eqref{5.13}.
That is, we wish to show that
\begin{equation}\label{5.14}
r^{1/2}\langle t-r\rangle |L^\nu Z^\alpha u'(t,x)|\le C\varepsilon
(1+t)^{3/20}\log(2+t)
\end{equation}
and
\begin{equation}\label{5.15}
\|\langle t+r\rangle L^\nu Z^\alpha \Box u(t,\cd)\|_2\le C\varepsilon
(1+t)^{3/20}\log(2+t)
\end{equation}
for $\nu \le 2$ and $|\alpha|+\nu \le 63$.
Notice that the first follows from the second by \eqref{4.7},
\eqref{5.4}, \eqref{5.6}, \eqref{5.10}, and \eqref{5.12}.  For \eqref{5.15}, we
expand $\Box u$
according to \eqref{1.7} to see that the left side is dominated by
$$\Bigl\|\Bigl(\langle t+r\rangle \sum_{\substack{|\alpha|+\mu\le
32\\\mu\le 2}} |L^\mu Z^\alpha u'(t,\cd)|\Bigr)
\sum_{\substack{|\alpha|+\mu\le 64\\ \mu\le 2}}|L^\mu Z^\alpha
u'(t,\cd)|\Bigr\|_2.$$
By \eqref{5.10} and \eqref{5.12}, this is easily seen to be bounded by
the right side of \eqref{5.15} as desired.

Since $|\Box (w-v)|\le C|\Box u|$, it is clear that the same proof
also yields
\begin{equation}\label{5.16}
r^{1/2}\langle t-r\rangle |L^\nu Z^\alpha (w-v)'(t,x)|\le C\varepsilon
(1+t)^{3/20}\log(2+t).
\end{equation}

Let's begin with (i.).  Since $v$ satisfies the better bound
\eqref{5.6}, it suffices to show
\begin{equation}\label{5.17}
\sum_{\substack{|\alpha|+\nu\le 52\\\nu\le 2}} \|L^\nu Z^\alpha
(w-v)'(t,\cdot)\|^2_2 \le C\varepsilon^3.
\end{equation}

By the standard energy integral method (see, e.g., Sogge \cite{S},
p.12), we have that the left side of \eqref{5.17} is bounded by
\begin{multline*}
C\sum_{\substack{|\alpha|+\nu\le
52\\\nu\le 2}}\int_0^t \int_\ext \Bigl|\Bigl\langle  
\partial_0 L^\nu Z^\alpha  (w-v),\Box L^\nu Z^\alpha(w-v)\Bigr\rangle\Bigr| \:dy\:ds \\
+ C \sum_{\substack{|\alpha|+\nu\le 52\\\nu\le 2}}\Bigl|\int_0^t \sum_{a=1}^3
\int_{\partial\mathcal{K}}  \partial_0 L^\nu Z^\alpha (w-v)
\partial_a L^\nu Z^\alpha (w-v) \:n_a \:d\sigma \:ds \Bigr|
\end{multline*}
where $n=(n_1,n_2,n_3)$ is the outward normal at a given point on
$\partial\mathcal{K}$ and $\langle \cd,\cd\rangle$ is the standard
Euclidean inner product on $\R^D$. Since $\mathcal{K}\subset
\{|x|<1\}$ and since $|L^\nu Z^\alpha (w-v)'(t,x)|\le
C\sum_{|\beta|\le |\alpha|, \mu\le \nu}|L^\mu \partial^\beta (w-v)'(t,x)|$ for
$x\in\partial\mathcal{K}$, we have that the last term is bounded by
$$C\int_0^t \int_{\{x\in \ext : |x|<1\}}\sum_{\substack{|\alpha|+\nu\le
53\\\nu\le
2}}|L^\nu \partial^\alpha (w-v)'(s,y)|^2\:dy\:ds.$$  Since we also
have that $[\Box, L]=2\Box$ and $[\Box, Z]=0$ and that $\Box
(w-v)=(1-\eta)\Box u = (1-\eta)Q(du,d^2u)$, we see that the left
side of \eqref{5.17} is thus controlled by
\begin{multline*}
C\int_0^t \int_\ext \Bigl|\Bigl\langle \sum_{\substack{|\alpha|+\nu\le
52\\\nu\le 2}} \partial_0 L^\nu Z^\alpha  (w-v),
\sum_{\substack{|\alpha|+\nu\le 52\\\nu\le 2}} L^\nu Z^\alpha
Q(du,d^2u)\Bigr\rangle\Bigr| \:dy\:ds
\\+ C\int_0^t \int_{\{x\in \ext : |x|<1\}}\sum_{\substack{|\alpha|+\nu\le
53\\\nu\le
2}}|L^\nu \partial^\alpha (w-v)'(s,y)|^2\:dy\:ds.
\end{multline*}
By \eqref{1.7}, this is dominated by:
\begin{multline}\label{5.18}
C\int_0^t \int_\ext \Bigl|\sum_{K=1}^D \sum_{\substack{|\alpha|+\nu\le
52\\\nu\le 2}}
\partial_0 L^\nu Z^\alpha (w-v)^K \sum_{0\le j,k,l\le 3}
B^{KK,jk}_{K,l} \sum_{\substack{|\alpha|+\nu\le 52\\\nu\le 2}}
\partial_l L^\nu Z^\alpha  u^K \\ \times\sum_{\substack{|\alpha|+\nu\le
52\\\nu\le 2}}
\partial_j\partial_k L^\nu Z^\alpha  u^K\Bigr|\:dy\:ds
\\+C \int_0^t \int_\ext \Bigl|\sum_{K=1}^D  \sum_{\substack{|\alpha|+\nu\le
52\\\nu\le 2}}
\partial_0 L^\nu Z^\alpha (w-v)^K \sum_{0\le j,k\le 3}
A^{K,jk}_{KK}  \sum_{\substack{|\alpha|+\nu\le 52\\\nu\le 2}}
\partial_j L^\nu Z^\alpha  u^K \\\times \sum_{\substack{|\alpha|+\nu\le
52\\\nu\le 2}}
\partial_k L^\nu Z^\alpha  u^K\Bigr|\:dy\:ds
\\+C \int_0^t \int_\ext \sum_{\substack{1\le I,J,K\le D\\ (I,K)\ne
(K,J)}} \sum_{\substack{|\alpha|+\nu\le 52\\\nu\le 2}} |L^\nu
Z^\alpha  \partial(w-v)^K| \sum_{\substack{|\alpha|+\nu\le
52\\\nu\le 2}} |L^\nu Z^\alpha
\partial u^I| \\\times\sum_{\substack{|\alpha|+\nu\le 53\\\nu\le 2}}
|L^\nu Z^\alpha
\partial u^J|\:dy\:ds
\\+ C\int_0^t \int_{\{x\in \ext : |x|<1\}}\sum_{\substack{|\alpha|+\nu\le
53\\\nu\le
2}}|L^\mu \partial^\alpha (w-v)'(s,y)|^2\:dy\:ds.
\end{multline}

The first two terms in \eqref{5.18} satisfy the bounds of Lemma
\ref{lemma4.1}. The third term involves interactions between waves
of different speeds.

When dealing with the first three terms of \eqref{5.18}, depending
on the linear estimates we shall employ, at times we shall use
certain $L^2$ and $L^\infty$ bounds for $u$ while at other times,
we shall use them for $w-v$. Since $u=(w-v)+v+u_0$ and $u_0$, $v$
satisfy the bounds \eqref{5.4},\eqref{5.6} respectively, it will
always be the case that bounds for $w-v$ will imply those for $u$
and vice versa.

Let us first handle the null terms.  By \eqref{4.1} and
\eqref{4.2}, the first two terms in \eqref{5.18} are controlled by
\begin{multline}\label{5.19}
C\int_0^t \int_\ext \sum_{\substack{|\alpha|+\mu\le 54\\\mu\le 3}} |L^\mu
Z^\alpha u| \sum_{\substack{|\alpha|+\mu\le 54\\\mu\le 3}} |L^\mu
Z^\alpha u'|\sum_{\substack{|\alpha|+\mu\le 52\\\mu\le 2}} |L^\mu
Z^\alpha (w-v)'| \:\frac{dy\:ds}{|y|}
\\+C \int_0^t \int_\ext \sum_{J=1}^D \frac{\langle c_Js-r\rangle}{\langle
s+r\rangle}
\sum_{\substack{|\alpha|+\mu\le 52\\\mu\le 2}} |L^\mu Z^\alpha
\partial (w-v)|\Bigl(\sum_{\substack{|\alpha|+\mu\le 53\\\mu\le 2}} |L^\mu
Z^\alpha
\partial u|\Bigr)^2\:dy\:ds
\end{multline}

To handle the contribution of the first term of
\eqref{5.19}, notice that by \eqref{5.4},\eqref{5.6}, and
\eqref{5.10} we have
$$\sum_{\substack{|\alpha|+\mu \le 54\\ \mu\le 3}} |L^\mu Z^\alpha
u(s,y)|\le C\varepsilon \langle s+|y|\rangle^{-9/10}\log(2+s),$$ 
which means
that the first term of \eqref{5.19} has a
contribution to \eqref{5.18} which is dominated by
\begin{multline*}
C\varepsilon \int_0^t  \frac{\log(2+s)}{\langle s\rangle^{9/10}}
\sum_{\substack{|\alpha|+\mu\le 54\\\mu\le
3}}\|\langle y\rangle^{-1/2} L^\mu
Z^\alpha u'(s,y)\|_2\\\times  \sum_{\substack{|\alpha|+\mu\le 52\\\mu\le
2}} \|\langle y\rangle^{-1/2} L^\mu Z^\alpha (w-v)'(s,y)\|_2\:ds
\end{multline*}
by the Schwarz inequality.  Thus, if we again apply the Schwarz inequality
and \eqref{5.13}, we see that this contribution is
$O(\varepsilon^3)$.

We now want to show that the second term of
\eqref{5.19} satisfies a similar bound.  If we apply \eqref{5.16},
we see that the second term of \eqref{5.19} is
controlled by
\begin{multline}\label{5.20}
C\varepsilon \int_0^t (1+s)^{3/20}\log(2+s) \int_\ext
\frac{1}{r^{1/2}\langle s+r\rangle} \sum_{\substack{|\alpha|+\mu\le
53\\\mu\le 2}}|L^\mu Z^\alpha \partial u|^2 \:dy\:ds
\\ \le C\varepsilon \int_0^t \frac{\log(2+s)}{(1+s)^{27/20}}
\sum_{\substack{|\alpha|+\mu\le 53 \\\mu\le 2}}\|L^\mu Z^\alpha
u'(s,\cd)\|^2_{L^2(|y|>s/2)} \:ds \\
+C\varepsilon \int_0^t (1+s)^{3/20}\log(2+s) \int_{\ext,\  |y|\le
s/2}\frac{1}{r^{1/2}\langle s+r\rangle}
\sum_{\substack{|\alpha|+\mu\le 53\\\mu\le 2}}|L^\mu Z^\alpha
u'(s,y)|^2\:dy\:ds.
\end{multline}
The first term on the right of \eqref{5.20} is clearly
$O(\varepsilon^3)$ by \eqref{5.12}.  For the second term on the right
of \eqref{5.20}, we apply \eqref{5.14} to control it as follows
$$C\varepsilon^2 \int_0^t \frac{1}{\langle s\rangle^{(6/5)-2\delta}}\int_{\ext,  \ |y|\le s/2}
\frac{1}{r^{(3/2)+\delta}} \sum_{\substack{|\alpha|+\mu\le 53\\\mu\le
2}}|L^\mu Z^\alpha u'(s,y)|\:dy\:ds.$$
Thus, if $\delta$ is sufficiently small, the Schwarz inequality and
\eqref{5.12} show that this term is also $O(\varepsilon^3)$.  This
concludes the proof that the contribution of the null forms enjoys an
$O(\varepsilon^3)$ bound.

We now wish to show that the multi-speed terms
\begin{equation}\label{5.21}
\int_0^t\int_\ext \sum_{\substack{|\alpha|+\mu\le 52\\\mu\le 2}} |\partial L^\mu
Z^\alpha (w-v)^K|
\sum_{\substack{|\alpha|+\mu\le 52\\\mu\le 2}}
|\partial L^\mu Z^\alpha u^I| \sum_{\substack{|\alpha|+\mu\le
53\\\mu\le 2}} |\partial L^\mu Z^\alpha u^J|\:dy\:ds
\end{equation}
 with $(I,K)\neq (K,J)$
have the same contribution to \eqref{5.18}.  For simplicity,
let us assume that $I\ne K$, $I=J$.  A symmetric argument will yield the same
bound for the remaining cases.  If we set $\delta < |c_I-c_K|/2$, it
follows that $\{|y|\in [(1-\delta)c_Is,(1+\delta)c_Is]\}\cap \{|y|\in
[(1-\delta)c_Ks,(1+\delta) c_Ks]\}=\emptyset$.  Thus, it will suffice
to show the bound when the spatial integral is taken over the complements each
of
these sets separately.  We will show the bound over  $\{|y|
\not\in [(1-\delta)c_Ks,(1+\delta)c_Ks]\}$.  The same argument will
symmetrically yield the bound over the other set.

If we apply \eqref{5.16}, we see that over $\{|y|
\not\in [(1-\delta)c_Ks,(1+\delta)c_Ks]\}$ \eqref{5.21} is
bounded by
\begin{equation}\label{5.22}
\varepsilon\int_0^t \frac{\log(2+s)}{\langle s \rangle^{17/20}} \int_{\ext, \ |y|\not\in [(1-\delta)c_Ks,
(1+\delta)c_Ks]}
\frac{1}{ r^{1/2}}
\sum_{\substack{|\alpha|+\mu\le 53\\\mu\le 2}}
|\partial L^\mu Z^\alpha u^I|^2 \:dy\:ds.
\end{equation}
Arguing as above, it is easy to see that these multiple speed terms
are also $O(\varepsilon^3)$.

Finally, we need to show that the last term in \eqref{5.18} enjoys
an $O(\varepsilon^4)$ contribution.  This is clear, however, since
this term is bounded by
$$\int_0^t \sum_{\substack{|\alpha|+\mu\le 53\\\mu\le 2}}\|L^\mu
\partial^\alpha (w-v)'(s,\cd)\|^2_\infty \:ds.$$
An application of \eqref{5.10} yields the desired bounds and
completes the proof of (i.).


We are now ready to prove (ii).  That is, we want to show
that we can prove \eqref{5.9} with $A_1$ replaced by $A_1/2$.  In view
of the bounds \eqref{5.6}, we see that it suffices to prove
\begin{equation}\label{5.23}
(1+t+r)\sum_{|\alpha|\le 40}|Z^\alpha (w-v)'(t,x)|\le C\varepsilon^{3/2}.
\end{equation}

The estimate is straightforward when $|x|>t/10$.  For then if we use
Lemma \ref{sobolev} and \eqref{5.17}, we get
\begin{equation}\begin{split}\label{5.24}
(1+t+|x|)\sum_{|\alpha|\le 40}|Z^\alpha (w-v)'(t,x)|&\le
C\sum_{|\alpha|\le 42}\|Z^\alpha (w-v)'(t,\cd)\|_2\\
&\le C\varepsilon^{3/2}, \quad |x|>t/10.
\end{split}\end{equation}

On account of this we only need to estimate the left side of
\eqref{5.23} when $|x|<t/10$.  Notice that $\Box
(w-v)=(1-\eta)Q(du,d^2u)$ vanishes when $|x|>10t$.  Thus, we can
apply \eqref{3.3} to conclude that when $|x|<t/10$, the left
side of \eqref{5.23} is dominated by
\begin{eqnarray*}
&& \sum_{\substack{|\beta|+\mu\le 43 \\ \mu\le 1}}\int_{t/100}^t \int |L^\mu Z^\beta
\partial[(1-\eta)Q(du,d^2u)]|\:\frac{dy\:ds}{|y|}\\
&& +C\sup_{0\le s\le t} (1+s)\sum_{|\beta|\le 44}\|Z^\beta
[(1-\eta)Q(du,d^2u)](s,\cd)\|_\infty\\
&& +C\sup_{0\le s\le t}(1+s)\sum_{|\beta|\le 47}\int_0^s
\int_{\substack{||y|-(s-\tau)| < 10\\ |y|\le
(600+\tau)/2}}|Z^\beta [(1-\eta)Q(du,d^2u)](\tau,y)|\:\frac{dy\:d\tau}{|y|}
\\
&& +C \sup_{0\le s\le t} \sum_{\substack{|\beta|+\mu\le 48\\\mu\le
1}}\frac{1}{1+s} \int_0^s\int_{|y|\ge (1+\tau)/10}|L^\mu Z^\beta
[(1-\eta)Q(du,d^2u)](\tau,y)| \:dy\:d\tau\\
&& =I+II+III+IV.
\end{eqnarray*}

Terms $II$ and $IV$ are the easiest to handle.  Since $u=(w-v)+v+u_0$,
by using \eqref{5.4},\eqref{5.6}, and \eqref{5.10}, one finds that $II$
is $O(\varepsilon^2)$.  Additionally, since \eqref{5.4}, \eqref{5.6}, and
\eqref{5.8} yield
$$\sum_{\substack{|\beta|+\mu\le 49\\\mu\le 1}}\|L^\mu Z^\beta
u'(\tau,\cd)\|_2\le C\varepsilon,$$
we can conclude that $IV$ is also $O(\varepsilon^2)$.

Similar considerations imply that
\begin{equation}\label{5.25}
I\le \int_{t/100}^t \int_{|y|<s/2} \sum_{\substack{|\beta|+\mu\le 43\\ \mu\le
1}}|L^\mu Z^\beta \partial Q(du,d^2u)(s,y)|\:\frac{dy\:ds}{|y|}+ C\varepsilon^2,
\end{equation}
since $|y|^{-1}=O(1/t)$ when $|y|>s/2$ and $t/100<s<t$ and since
$\partial \eta = O(1/t)$ when $t/100<s<t$.  The first term on the
right side of \eqref{5.25} is dominated by
\begin{equation}\label{5.26}
\int_{t/100}^t \int_{|y|<s/2}\sum_{\substack{|\beta|+\mu\le 44 \\ \mu\le 1}}|r^{1/2} \langle
s-r\rangle L^\mu Z^\beta u'| \sum_{\substack{|\beta|+\mu\le 44 \\ \mu \le 1}}
|r \langle s-r\rangle L^\mu Z^\beta u''|\:\frac{dy\:ds}{|y|^3
s^{3/2}}.\end{equation}
If we apply Lemma \ref{sobolev} and \eqref{4.3}, we see that
\begin{equation}\label{5.27}
\sum_{\substack{|\beta|\le 44 \\\mu\le 1}} |r\langle s-r\rangle L^\mu Z^\beta
u''(s,y)|\le C\varepsilon (1+s)^{3/20}\log(2+s)
\end{equation}
by \eqref{5.4}, \eqref{5.6}, \eqref{5.8}, and \eqref{5.15}.  Thus,
it follows that \eqref{5.26}, and hence $I$, is also $O(\varepsilon^2)$ using
\eqref{5.14} and \eqref{5.27}.

It remains to estimate $III$.  If we use \eqref{4.7}, we
conclude that on the region of integration
\begin{equation}\begin{split}\label{5.28}
&\sum_{|\beta|\le 48}|Z^\beta u'(\tau,y)| \\
&\le \frac{C}{r^{1/2}\tau}\Bigl[\sum_{\substack{|\beta|+\mu\le
50\\\mu\le 1}}\|L^\mu Z^\beta u'(\tau,\cd)\|_2 + \sum_{|\beta|\le
49}\|\langle \tau+r\rangle Z^\beta \Box
u(\tau,\cd)\|_2+(1+\tau)\|u'(\tau,\cd)\|_\infty\Bigr]\\
&\le \frac{C}{r^{1/2}\tau}\Bigl[\varepsilon + \sum_{|\beta|\le 49}
\|\langle \tau+r\rangle Z^\beta \Box u(\tau,\cd)\|_2\Bigr],
\end{split}\end{equation}
using \eqref{5.4},\eqref{5.6},\eqref{5.8}, and \eqref{5.9} in the last step.
If $|\beta|\le 49$
$$\langle \tau+r\rangle |Z^\beta \Box u(\tau,y)|\le \sum_{|\gamma|\le
50}|Z^\gamma u'(\tau,y)| \times \Bigl(\langle \tau+r\rangle
\sum_{|\gamma|\le 25} |Z^\gamma u'(\tau, y)|\Bigr),$$
which by the low energy estimate \eqref{5.8} and the low dispersive
estimate \eqref{5.9} gives
$$\|\langle \tau+r\rangle Z^\beta \Box u(\tau,\cd)\|_2 \le
C\varepsilon \sup_y \langle \tau+r\rangle \sum_{|\gamma|\le
25}|Z^\gamma u'(\tau, y)| \le C\varepsilon^2.$$
Combining this with \eqref{5.28} and recalling that
$(1-\eta(\tau,y))=0$ for $|y|>10\tau$, we get
\begin{equation}\label{5.29}
III\le C\varepsilon^2 \sup_{0\le s\le t}(1+s)\int_{s/100}^s
\int_{|r-(s-\tau)|<10} \Bigl(\frac{1}{r^{1/2}\tau}\Bigr)^2
r\:dr\:d\tau \le C\varepsilon^2
\end{equation}
which completes the proof of (ii.).

To complete the proof of Theorem \ref{theorem1.1}, we need to show
how \eqref{5.8},\eqref{5.9} imply \eqref{5.10}-\eqref{5.13}.

Since \eqref{5.9} has been established, the remainder of the
argument follows nearly verbatim from the arguments of \cite{MS}.
For completeness, we will sketch the argument here.  We begin by
using the above facts to prove \eqref{5.11}.  With notation as in \S
1-2, $\Box_\gamma u=B(du)$ with
$$\gamma^{IJ,jk}=-\sum_{\substack{0\le l\le 3\\1\le K\le D}}
B^{IJ,jk}_{K.l}\partial_l u^K.$$  By \eqref{5.9}, we have
$$\|\gamma'(s,\cd)\|_\infty \le \frac{C\varepsilon}{1+s}.$$
Let us first show
the estimates for the energy of $\partial_t^j u$ for $j\le M\le
100$. We shall use induction on $M$.

We first notice that by \eqref{2.5} and \eqref{5.9} we have
\begin{equation}\label{5.30}
\partial_t E^{1/2}_M(u)(t)\le C\sum_{j\le M}\|\Box_\gamma
\partial^j_t u(t,\cd)\|_2 +
\frac{C\varepsilon}{1+t}E^{1/2}_M(u)(t).
\end{equation}
Note that for $M=1,2,\dots$
\begin{align*}
\sum_{j\le M}|\Box_\gamma \partial^j_t u| &\le C\Bigl(\sum_{j\le
M} |\partial^j_t u'|+\sum_{j\le M-1}|\partial_t^j \partial^2
u|\Bigr)\sum_{|\alpha|\le 40}|\partial^\alpha  u'| \\
&\qquad\qquad\qquad\qquad + C\sum_{|\alpha|\le M-40}
|\partial^\alpha u'|\sum_{41\le|\alpha|\le M/2}|\partial^\alpha
u'|\\
&\le \frac{C\varepsilon}{1+t}\Bigl(\sum_{j\le M}|\partial^j_t
u'|+\sum_{j\le M-1}|\partial^j_t\partial^2 u|\Bigr) + C
\sum_{|\alpha|\le M-40}|\partial^\alpha u'|\sum_{41\le|\alpha|\le
M/2}|\partial^\alpha u'|
\end{align*}
by \eqref{5.9} and \eqref{5.4}.  Also, if we use elliptic
regularity and repeat this argument, we get
\begin{align*}
\sum_{j\le M-1}\|\partial_t^j \partial^2 u(t,\cd)\|_2 &\le C
\sum_{j\le M}\|\partial_t^j u'(t,\cd)\|_2 + C\sum_{j\le
M-1}\|\partial_t^j \Box u(t,\cd)\|_2\\
&\le C\sum_{j\le M}\|\partial_t^j u'(t,\cd)\|_2 +
\frac{C\varepsilon}{1+t}\sum_{j\le M-1}\|\partial^j_t\partial^2
u(t,\cd)\|_2\\
&\quad\quad\quad + C \sum_{|\alpha|\le M-41, |\beta|\le M/2}
\|\partial^\alpha u'(t,\cd)\partial^\beta u'(t,\cd)\|_2.
\end{align*}
If $\varepsilon$ is small, we can absorb the second to last term
into the left side of the preceding inequality.  Therefore, if we
combine the last two inequalities, we conclude that
\begin{multline*}
\sum_{j\le M} \|\Box_\gamma \partial_t^j u(t,\cd)\|_2\le
\frac{C\varepsilon}{1+t}\sum_{j\le M}\|\partial^j_t u'(t,\cd)\|_2
\\+ C \sum_{|\alpha|\le M-40, |\beta|\le M/2}\|\partial^\alpha
u'(t,\cd)\partial^\beta u'(t,\cd)\|_2.
\end{multline*}
If we combine this with \eqref{5.30} we get that for small
$\varepsilon>0$
\begin{equation}\label{5.31}
\partial_t E_M^{1/2}(u)(t)\le
\frac{C\varepsilon}{1+t}E^{1/2}_M(u)(t)+C\sum_{|\alpha|\le M-40,
|\beta|\le M/2}\|\partial^\alpha u'(t,\cd)\partial^\beta
u'(t,\cd)\|_2,
\end{equation}
since when $\varepsilon$ is small, $\frac{1}{2}E^{1/2}_M(u)(t)\le
\sum_{j\le M}\|\partial^j_t u'(t,\cd)\|_2\le 2E^{1/2}_M(u)(t)$.

For $M\le 52$, the energy estimate \eqref{5.11} follows from
\eqref{5.8}.  When $M>52$ we have to deal with the last term in
\eqref{5.31}.  To do this we first note that by Lemma
\ref{sobolev} we have
$$\sum_{|\alpha|\le M-40, |\beta|\le
M/2}\|\partial^{\alpha}u'(t,\cd) \partial^\beta u'(t,\cd)\|_2\le C
\sum_{|\gamma|\le \max(M-38,2+M/2)}\|\langle
x\rangle^{-1/2}Z^\gamma u'(t,\cd)\|^2_2,$$ which means that for
$40\le M\le 100$, \eqref{5.31}, \eqref{5.1}, and Gronwall's
inequality yield
\begin{equation}\label{5.32}
E_M^{1/2}(u)(t)\le C (1+t)^{C\varepsilon}\Bigl[\varepsilon + 
\sum_{|\alpha|\le \max(M-38,2+M/2)}\|\langle
x\rangle^{-1/2}Z^\alpha u'\|^2_{L^2(S_t)}\Bigr],
\end{equation}
if, as before, $S_t=[0,t]\times \ext$.

If we use \eqref{5.8} and \eqref{5.32} along with a simple
induction argument we conclude that we would have the desired
bounds
\begin{equation}\label{5.33}
E_{100}^{1/2}(u)(t)\le C\varepsilon (1+t)^{C\varepsilon +
\sigma}\end{equation}
 for arbitrarily small $\sigma>0$ if we apply
the following lemma.
\begin{lemma}\label{lemma5.2}
Under the above assumptions, if $M\le
100-8\mu$, $\mu \le 4$, and
\begin{multline}\label{5.34}
\sum_{\substack{|\alpha|+\nu\le M\\\nu\le \mu}}\|L^\nu \partial^\alpha u'(t,\cd)\|_2 +
\sum_{\substack{|\alpha|+\nu\le M-3\\\nu\le \mu}}\|\langle x\rangle^{-1/2} L^\nu
\partial^\alpha u'\|_{L^2(S_t)}\\+\sum_{\substack{|\alpha|+\nu\le M-4\\\nu\le\mu}}\|L^\nu
Z^\alpha u'(t,\cd)\|_2
+\sum_{\substack{|\alpha|+\nu\le M-6\\\nu\le\mu}}\|\langle x\rangle^{-1/2} L^\nu Z^\alpha
u'\|_{L^2(S_t)}\le C\varepsilon (1+t)^{C\varepsilon+\sigma},
\end{multline}
with $\sigma>0$, then there is a constant $C'$ so that
\begin{multline}\label{5.35}
\sum_{\substack{|\alpha|+\nu\le M-2\\\nu\le\mu}}\|\langle x\rangle^{-1/2}L^\nu
\partial^\alpha u'\|_{L^2(S_t)}+\sum_{\substack{|\alpha|+\nu\le M-3\\\nu\le\mu}}\|L^\nu
Z^\alpha u'(t,\cd)\|_2 \\
+\sum_{\substack{|\alpha|+\nu\le M-5\\\nu\le\mu}}\|\langle x\rangle^{-1/2}L^\nu Z^\alpha
u'\|_{L^2(S_t)}\le C'\varepsilon (1+t)^{C'\varepsilon + C'\sigma}.
\end{multline}
\end{lemma}
The proof of this lemma can be found in \cite{MS}.

By elliptic regularity and \eqref{5.33}, we get \eqref{5.11}.
Also, from Lemma \ref{lemma5.2}, we get
\begin{multline}\label{5.36}
\sum_{|\alpha|\le 98}\|\langle x\rangle^{-1/2}\partial^\alpha
w'\|_{L^2(S_t)}+\sum_{|\alpha|\le 97}\|Z^\alpha w'(t,\cd)\|_2\\
+\sum_{|\alpha|\le 95}\|\langle x\rangle^{-1/2} Z^\alpha
w'\|_{L^2(S_t)}\le C'\varepsilon (1+t)^{C'\varepsilon+C'\sigma},
\end{multline}
since the same sort of bounds hold when $w$ is replaced by $u$.

Here and in what follows $\sigma$ denotes a small constant that
must be taken to be larger and larger at each occurrence.  Note
that in terms of the number of $Z$ derivatives \eqref{5.35} is
considerably stronger than the variants of \eqref{5.12} and
\eqref{5.13} where one just takes the terms with $\nu=0$.  This is
because just as there is a loss of derivatives 
in going from \eqref{5.11} to \eqref{5.36}, there will also be a loss of derivatives in
going from $L^2$ bounds for terms of the form $L^\nu Z^\alpha u'$
to those of the form $L^{\nu+1} Z^\alpha u'$.

The proof of the estimates involving powers of $L$ is a bit more
complicated, but still follows the strategy above.  First we will
estimate $L^\nu \partial^\alpha u'$ in $L^2$ when $\alpha$ is
small using \eqref{5.9}.  Then we shall estimate the remaining
parts of \eqref{5.12} and \eqref{5.13} using Lemma \ref{lemma5.2}.

The main part of the next step is to show that
\begin{equation}\label{5.37}
\sum_{\substack{|\alpha|+\mu\le 92\\ \mu\le 1}} \|L^\mu
\partial^\alpha u'(t,\cd)\|_2\le C\varepsilon (1+t)^{C\varepsilon
+ \sigma}.
\end{equation}

For this we shall want to use \eqref{2.8}.  We must first
establish appropriate versions of \eqref{2.7} for $N_0+\nu_0\le
92, \nu_0=1$.  For this we note that for $M\le 92$
\begin{multline*}
\sum_{\substack{j+\mu\le M\\\mu \le 1}}\Bigl(|\tilde{L}^\mu \partial_t^j \Box_\gamma
 u| + |[\tilde{L}^\mu \partial_t^j, \Box -\Box_\gamma]u|\Bigr)\\
\le C\Bigl(\sum_{j\le M-1}|\tilde{L}\partial_t^j
\partial u|+\sum_{j\le M-2}|\tilde{L}\partial_t^j
\partial^2 u|\Bigr)\sum_{|\alpha|\le 40}|\partial^\alpha u'| \\+
C\sum_{|\alpha|\le M-41}|L\partial^\alpha u'|\sum_{|\alpha|\le
M}|\partial^\alpha u'| + C\sum_{|\alpha|\le M}|\partial^\alpha
u'|\sum_{|\alpha|\le \max(M/2,M-40)}|\partial^\alpha u'|.
\end{multline*}
By this, \eqref{5.9}, and elliptic regularity, we get that for
$M\le 92$
\begin{multline*}
\sum_{\substack{j+\mu\le M\\ \mu\le 1}}\Bigl(\|\tilde{L}^\mu \partial_t^j
\Box_\gamma u(t,\cd)\|_2 + \|[\tilde{L}^\mu \partial_t^j, \Box-\Box_\gamma]u(t,\cd)\|_2\Bigr)
\le
\frac{C\varepsilon}{1+t}\sum_{\substack{j+\mu\le M\\\mu\le
1}}\|\tilde{L}^\mu \partial_t^j u'(t,\cd)\|_2
\\+C \sum_{|\alpha|\le M-41}\|\langle x\rangle^{-1/2}
L\partial^\alpha u'(t,\cd)\|_2\sum_{|\alpha|\le 94}\|\langle
x\rangle^{-1/2}Z^\alpha u'(t,\cd)\|_2 \\
+C\sum_{|\alpha|\le \max (M, 2+M/2)}\|\langle x\rangle^{-1/2}
Z^\alpha u'(t,\cd)\|_2^2.
\end{multline*}
Based on this if $\varepsilon$ is small then \eqref{2.7} holds
with $\delta=C\varepsilon$ and
$$H_{1,M-1}(t)=C\sum_{|\alpha|\le M-41}\|\langle
x\rangle^{-1/2}L\partial^\alpha u'(t,\cd)\|^2_2+C\sum_{|\alpha|\le
94}\|\langle x\rangle^{-1/2}Z^\alpha u'(t,\cd)\|^2_2.$$ 
Therefore
since the conditions on the data give $\int e_0(\tilde{L}^\nu
\partial^j_t u)(0,x)\:dx \le C\varepsilon^2$ if $\nu+j\le 100$ it
follows from \eqref{2.8} and \eqref{5.36} that for $M\le 92$
\begin{multline}\label{5.38}
\sum_{\substack{|\alpha|+\mu\le M\\\mu\le 1}}\|L^\mu
\partial^\alpha u'(t,\cd)\|_2 \le C\varepsilon
(1+t)^{C\varepsilon+\sigma}+C(1+t)^{C\varepsilon}\sum_{|\alpha|\le
M-41}\|\langle x\rangle^{-1/2}L\partial^\alpha u'\|^2_{L^2(S_t)}
\\+ C(1+t)^{C\varepsilon}\int_0^t \sum_{|\alpha|\le
M+1}\|\partial^\alpha u'(s,\cd)\|_{L^2(|x|<1)}\:ds.
\end{multline}
If we apply \eqref{2.13} and \eqref{5.4} we get that the last
integral is dominated by $\varepsilon \log(2+t)$ plus
\begin{multline*}
\int_0^t \sum_{|\alpha|\le M+1}\|\partial^\alpha
w'(s,\cd)\|_{L^2(|x|<1)}\:ds
\\\le C \sum_{|\alpha|\le M+2}\int_0^t \Bigl(\int_0^s
\|\partial^\alpha \Box
w(\tau,\cd)\|_{L^2(||x|-(s-\tau)|<10)}\:d\tau\Bigr)\:ds.
\end{multline*}
By \eqref{5.4} if we replace $w$ by $u_0$ then the analog of the
last term is $O(\varepsilon \log(2+t))$.  We therefore conclude
that
\begin{multline*}
\sum_{|\alpha|\le M+1}\int_0^t \|\partial^\alpha
u'(s,\cd)\|_{L^2(|x|<1)}\:ds \le C\varepsilon \log(2+t) \\
+C\sum_{|\alpha|\le M+2}\int_0^t \Bigl(\int_0^s \|\partial^\alpha
\Box u(\tau,\cd)\|_{L^2(||x|-(s-\tau)|< 10)}\:d\tau\Bigr)\:ds.
\end{multline*}
Since
$$\sum_{|\alpha|\le M+2}|\partial^\alpha \Box u|\le
C\sum_{|\alpha|\le M+3}|\partial^\alpha u'|\sum_{|\alpha|\le
(M+3)/2}|\partial^\alpha u'|,$$ 
an application of Lemma
\ref{sobolev} yields
$$\sum_{|\alpha|\le M+2}\|\partial^\alpha \Box
u(\tau,\cd)\|_{L^2(||x|-(s-\tau)|<10)}\le C\sum_{|\alpha|\le
95}\|\langle x\rangle^{-1/2}Z^\alpha
u'\|^2_{L^2(||x|-(s-\tau)|<20)},$$ 
since $(3+M)/2\le 95$ if $M\le
92$.  Since the sets $\{(\tau,x):||x|-(j-\tau)|<20\}$,
$j=0,1,2,\dots$ have finite overlap, we conclude that for $M\le
92$
\begin{align*}
\sum_{|\alpha|\le M+1}\int_0^t
\|\partial^{\alpha}u'(s,\cd)\|_{L^2(|x|<1)}\:ds &\le C \varepsilon
\log(2+t)+C\sum_{|\alpha|\le 95}\|\langle x\rangle^{-1/2} Z^\alpha
u'\|^2_{L^2(S_t)}\\
&\le C\varepsilon (1+t)^{C\varepsilon+\sigma}.
\end{align*}
Therefore, by \eqref{5.38} we have that
\begin{equation*}
\sum_{\substack{|\alpha|+\mu\le M\\\mu\le 1}}\|L^\mu
\partial^\alpha u'(t,\cd)\|_2\le C\varepsilon
(1+t)^{C\varepsilon+\sigma}+C(1+t)^{C\varepsilon}\sum_{|\alpha|\le
M-41}\|\langle x\rangle^{-1/2}L\partial^\alpha u'\|^2_{L^2(S_t)}.
\end{equation*}
This gives the desired bounds when $M\le 40$.

If we now use induction and Lemma \ref{lemma5.2}, we get \eqref{5.37} as
well as
\begin{multline}\label{5.39}
\sum_{\substack{|\alpha|+\mu\le 90\\\mu\le 1}}\|\langle
x\rangle^{-1/2}L^\mu \partial^\alpha
u'\|_{L^2(S_t)}+\sum_{\substack{|\alpha|+\mu\le 89\\\mu\le
1}}\|L^\mu Z^\alpha u'(t,\cd)\|_2\\
+\sum_{\substack{|\alpha|+\mu\le 87\\\mu\le 1}}\|\langle
x\rangle^{-1/2}L^\mu Z^\alpha u'\|_{L^2(S_t)}\le C\varepsilon
(1+t)^{C\varepsilon + C\sigma}.
\end{multline}

If we repeat this argument we can estimate $L^2 Z^\alpha u'$, $L^3
Z^\alpha u'$, and $L^4 Z^\alpha u'$ for appropriate $Z^\alpha$.
Using \eqref{5.37}, \eqref{5.39}, and the last argument gives
\begin{multline*}
\sum_{\substack{|\alpha|+\mu\le 84\\\mu\le 2}}\|L^\mu
\partial^\alpha u'(t,\cd)\|_2
+\sum_{\substack{|\alpha|+\mu\le 82\\ \mu\le 2}}\|\langle
x\rangle^{-1/2}L^\mu \partial^\alpha u'\|_{L^2(S_t)}
+\sum_{\substack{|\alpha|+\mu\le 81\\
\mu\le 2}}\|L^\mu Z^\alpha u'(t,\cd)\|_2 \\
+\sum_{\substack{|\alpha|+\mu\le 79\\ \mu\le 2}}\|\langle
x\rangle^{-1/2}L^\mu Z^\alpha u'\|_{L^2(S_t)}\le C\varepsilon
(1+t)^{C\varepsilon + C\sigma}.
\end{multline*}
Then using the estimates for $L^\mu Z^\alpha u'$, $\mu\le 2$ we
can argue as above to get
\begin{multline*}
\sum_{\substack{|\alpha|+\mu\le 76\\\mu\le 3}}\|L^\mu
\partial^\alpha u'(t,\cd)\|_2
+\sum_{\substack{|\alpha|+\mu\le 74\\ \mu\le 2}}\|\langle
x\rangle^{-1/2}L^\mu \partial^\alpha u'\|_{L^2(S_t)}
+\sum_{\substack{|\alpha|+\mu\le 73\\
\mu\le 3}}\|L^\mu Z^\alpha u'(t,\cd)\|_2 \\
+\sum_{\substack{|\alpha|+\mu\le 71\\ \mu\le 3}}\|\langle
x\rangle^{-1/2}L^\mu Z^\alpha u'\|_{L^2(S_t)}\le C\varepsilon
(1+t)^{C\varepsilon + C\sigma}.
\end{multline*}
Similarly, using the estimates for $L^\mu Z^\alpha u'$ for $\mu\le
3$ we finally get
\begin{multline*}
\sum_{\substack{|\alpha|+\mu\le 68\\\mu\le 4}}\|L^\mu
\partial^\alpha u'(t,\cd)\|_2
+\sum_{\substack{|\alpha|+\mu\le 66\\ \mu\le 2}}\|\langle
x\rangle^{-1/2}L^\mu \partial^\alpha u'\|_{L^2(S_t)}
+\sum_{\substack{|\alpha|+\mu\le 65\\
\mu\le 4}}\|L^\mu Z^\alpha u'(t,\cd)\|_2 \\
+\sum_{\substack{|\alpha|+\mu\le 63\\ \mu\le 4}}\|\langle
x\rangle^{-1/2}L^\mu Z^\alpha u'\|_{L^2(S_t)}\le C\varepsilon
(1+t)^{C\varepsilon + C\sigma}.
\end{multline*}
If we combine this with our earlier bounds, we conclude that
\eqref{5.12} and \eqref{5.13} must be valid.

It remains to prove \eqref{5.10}.  This is straightforward.  If we
use Theorem \ref{theorem3.1} we find that its left side is
dominated by the square of that of \eqref{5.13}.  Hence
\eqref{5.13} implies \eqref{5.10} which completes the proof.


\end{document}